\documentclass[12pt,oneside,american]{amsart}
\usepackage[T1]{fontenc}
\usepackage[utf8]{inputenc}
\usepackage{xcolor}
\usepackage{mathrsfs}
\usepackage{amstext}
\usepackage{amsthm}
\usepackage{amssymb}
\PassOptionsToPackage{normalem}{ulem}
\usepackage{ulem}
\usepackage{microtype}

\makeatletter

\providecolor{lyxadded}{rgb}{0,0,1}
\providecolor{lyxdeleted}{rgb}{1,0,0}
\DeclareRobustCommand{\mklyxadded}[1]{\textcolor{lyxadded}\bgroup#1\egroup}
\DeclareRobustCommand{\mklyxdeleted}[1]{\textcolor{lyxdeleted}\bgroup\mklyxsout{#1}\egroup}
\DeclareRobustCommand{\mklyxsout}[1]{\ifx\\#1\else\sout{#1}\fi}

\numberwithin{equation}{section}
\numberwithin{figure}{section}
\theoremstyle{plain}
\newtheorem{thm}{\protect\theoremname}[section]
\theoremstyle{plain}
\newtheorem{prop}[thm]{\protect\propositionname}
\theoremstyle{definition}
\newtheorem{defn}[thm]{\protect\definitionname}
\theoremstyle{remark}
\newtheorem{rem}[thm]{\protect\remarkname}
\theoremstyle{definition}
\newtheorem{example}[thm]{\protect\examplename}
\theoremstyle{plain}
\newtheorem{lem}[thm]{\protect\lemmaname}
\theoremstyle{plain}
\newtheorem{cor}[thm]{\protect\corollaryname}

\usepackage{fullpage}

\makeatother

\usepackage{babel}
\providecommand{\corollaryname}{Corollary}
\providecommand{\definitionname}{Definition}
\providecommand{\examplename}{Example}
\providecommand{\lemmaname}{Lemma}
\providecommand{\propositionname}{Proposition}
\providecommand{\remarkname}{Remark}
\providecommand{\theoremname}{Theorem}

\begin{document}
\title{On the Isomorphism of Semicircular Algebras}
\address{Department of Mathematics, UCLA, Los Angeles, CA 90095}
\email{shlyakht@math.ucla.edu}
\date{\today}
\author{Dima Shlyakhtenko}
\begin{abstract}
Building on OpenAI's proof of isomorphism of free group factors \cite{OpenAIFreeGroupFactors},
we give a self-contained proof of a corresponding result for $C^{*}$-algebras
generated by free semicircular systems: the isomorphism class of a
$C^{*}$-algebra generated by a semicircular $n$-tuple does not depend
on $n\in\{2,3,\dots,+\infty\}$. Similarly, reduced free products
of $n$ copies of the Jiang--Su algebra are isomorphic for different
$n$ in this range. While reduced free group $C^{*}$-algebras of
different ranks are not isomorphic (by a result of Pimsner--Voiculescu),
we show that they have strong free absorption properties. At the heart
of the proof is an automorphism approximation property for inclusions
of $C^{*}$- and $W^{*}$-algebras; a relative version of that property
was used to derive free complementation results of \cite{ShlyakhtenkoFuchsian}.
\end{abstract}

\maketitle

\section{Introduction.}

Let $\mathbb{F}_{n}$ be the free group on $n$ generators. A ground-breaking
result of Pimsner and Voiculescu \cite{PimsnerVoiculescu} snowed
that the reduced $C^{*}$-algebras $C_{\textrm{red}}^{*}(\mathbb{F}_{n})$
and $C_{\textrm{red}}^{*}(\mathbb{F}_{m})$ are non-isomorphic, whenever
$n\neq m$, $n,m\in\{2,3,\dots,+\infty\}$ (see also \cite{GaoKunnawalkamElayavalli}).
As a $C^{*}$-algebra, $C_{\textrm{red}}^{*}(\mathbb{F}_{n})$ is
the $n$-fold reduced free product of algebras $C(\mathbb{T})$ of
continuous functions on the circle, taken with the Haar measure. 

Consider instead a free semicircular system $S_{1},S_{2},\dots,S_{n}$
\cite{VoiculescuSymmetries}; then $Q_{n}=C^{*}(S_{1},\dots,S_{n})$
is the $n$-fold free product of algebras $C[-2,2]$ of continuous
functions on an interval, endowed with the semicircular measure. It
is therefore natural to ask whether the isomorphism class of $Q_{n}$
depends on $n$. Since the von Neumann algebra $W^{*}(S_{1},\dots,S_{n})$
is the free group factor $L(\mathbb{F}_{n})$, the question of isomorphism
$Q_{n}\cong Q_{m}$ for $n\neq m$ is formally stronger than that
of $L(\mathbb{F}_{n})\cong L(\mathbb{F}_{m}).$ 

In this paper, we show that OpenAI's proof of isomorphism of free
group factors \cite{OpenAIFreeGroupFactors} can be modified to give
a proof of isomorphism of the semicircular $C^{*}$-algebras $Q_{n}$
for all $n\geq2$. Our proof is self-contained and does not rely on
the results of \cite{OpenAIFreeGroupFactors}. We show (see \S\ref{sec:free-product-consequences}):
\begin{itemize}
\item $C^{*}(S_{1},\dots,S_{n})\cong C^{*}(S_{1},S_{2})$, for any $n\in\{2,3,\dots,+\infty\}$;
this implies $L(\mathbb{F}_{n})\cong L(\mathbb{F}_{2})$.
\item $C_{\textrm{red}}^{*}(\mathbb{F}_{m})*_{\textrm{red}}C^{*}(S_{1},\dots,S_{n})\cong C_{\textrm{red}}^{*}(\mathbb{F}_{m})$,
for any $m\in\{2,3,\dots,+\infty\}$ and $n\in\{1,2,\dots,+\infty\}$
\item Let $\mathcal{Z}$ denote the Jiang--Su algebra taken with its canonical
trace. Then for any $n\in\{2,3,\dots,+\infty\}$ the $n$-fold reduced
free product $\mathcal{Z}^{*n}\cong\mathcal{Z}*_{\textrm{red}}\mathcal{Z}\cong C^{*}(S_{1},S_{2})$. 
\item $C_{\textrm{red}}^{*}(\mathbb{F}_{m})*_{\textrm{red}}\mathcal{Z}^{*n}\cong C_{\textrm{red}}^{*}(\mathbb{F}_{m})$,
for any $m\in\{2,3,\dots,+\infty\}$ and $n\in\{1,2,\dots,+\infty\}$
\end{itemize}
These results show that the free entropy dimension invariants (\cite{VoiculescuEntropyII,VoiculescuEntropyV,ShlyakhtenkoPartialConjugates,ConnesShlyakhtenko})
cannot be $C^{*}$-algebra invariants. Indeed, all of these invariants
give the value $n$ for the $n$-tuple $S_{1},\dots,S_{n}$. Likewise,
one cannot expect that numeric $L^{2}$-invariants, such as $L^{2}$-Betti
numbers for tracial $*$-algebras \cite{ConnesShlyakhtenko} or the
$L^{2}$ Euler characteristic, are invariants not just of von Neumann
algebras, but even of $C^{*}$-algebras. Note, however, that the question
of whether $L(\Gamma)\cong L(\Lambda)$ implies that the $L^{2}$-Betti
numbers $\beta_{1}^{(2)}(\Gamma)$ and $\beta_{1}^{(2)}(\Lambda)$
must be simultaneously zero or non-zero remains open (cf. \cite{ShlyakhtenkoSofic}).

The isomorphism result $L(\mathbb{F}_{t})\cong L(\mathbb{F}_{\infty})$
proved in \cite{OpenAIFreeGroupFactors} also implies that there exists
a continuum many of groups that are $W^{*}$-equivalent (i.e., their
von Neumann algebras are isomorphic) and yet not orbit-equivalent.
Indeed, the groups $G_{t}$ constructed in \cite{PopaShlyakhtenkoGt}
satisfy $L(G_{t})\cong L(\mathbb{F}_{t})\cong L(\mathbb{F}_{2})$,
and yet $\beta_{1}^{(2)}(G_{t})=t-1$, so $G_{t}\nsim_{OE}G_{s}$
if $t\neq s$, by results of Gaboriau \cite{GaboriauL2}. Note that
$\mathbb{F}_{2}$ and $\mathbb{F}_{\infty}$ are not even measure-equivalent,
since their $\ell^{2}$ Betti numbers are not proportional with a
finite constant \cite{GaboriauL2}. In fact, the groups $(G_{t}\times\mathbb{F}_{2})*\mathbb{F}_{2}$
satisfy $\beta_{1}^{(2)}((G_{t}\times\mathbb{F}_{2})*\mathbb{F}_{2})=2$,
$\beta_{2}^{(2)}((G_{t}\times\mathbb{F}_{2})*\mathbb{F}_{2})=(t-1)$.
These pairs are not proportional for different values of $t$, and
so $(G_{t}\times\mathbb{F}_{2})*\mathbb{F}_{2}$ and $(G_{s}\times\mathbb{F}_{2})*\mathbb{F}_{2}$
cannot be measure-equivalent by \cite{GaboriauL2}. This gives a continuum
of non-measure-equivalent groups .

Our reformulation of OpenAI's proof is centered on an abstract Automorphism
Approximation Property (AAP) for inclusions of von Neumann algebras
or $C^{*}$-algebras (in the latter case we call this property CAAP),
which we describe below. A relative version of AAP appeared as an
approximation criteria in \cite{ShlyakhtenkoFuchsian} and was used
to prove free complementation of the commutator MASA in $L(\mathbb{F}_{2})$. 

\subsection{The approximation properties AAP and CAAP}

We begin by describing two related properties of inclusions of von
Neumann algebras and $C^{*}$-algebras. Very roughly, the properties
say that the automorphism group of the larger algebra acts so wildly
as to be able to move an arbitrary element of the larger algebra arbitrarily
close to the subalgebra, while almost fixing an arbitrary finite subset
of the subalgebra. Such a property automatically implies that the
algebra and its subalgebra are isomorphic; moreover, the inclusion
is elementary (see below and \S \ref{sec:elementary-inclusions}).
\begin{prop}
\label{prop:aap}Let $Q\subset N\subset M$ be inclusions of tracial
von Neumann algebras, with $M$ having separable predual. The following
conditions are equivalent:
\begin{enumerate}
\item For every finite subset $F\subset N$, every $z\in M$, and every
$\eta,\epsilon>0$, there are a trace-preserving automorphism $\alpha$
of $M$ fixing $Q$ pointwise and $y\in N$ such that 
\[
\max_{a\in F}\|\alpha(a)-a\|_{2}<\eta,\qquad\|z-\alpha(y)\|_{2}<\epsilon.
\]
\item For every finite subset $F\subset N$ and every $\epsilon>0$, there
is a surjective trace-preserving isomorphism $\Phi:N\to M$ fixing
$Q$ pointwise such that 
\[
\max_{a\in F}\|\Phi(a)-a\|_{2}<\epsilon.
\]
\item There is a normal trace-preserving endomorphism $\rho:M\to M$ with
$\rho(M)=N$ which is a pointwise $L^{2}$ limit of a sequence of
trace-preserving automorphisms of $M$ fixing $Q$ pointwise.
\end{enumerate}
Moreover, fix a unital von Neumann subalgebra $E\subset N$. Conditions
(1) and (2) remain equivalent if, for every additional tolerance $\kappa>0$,
the map $\alpha$, respectively $\Phi$, is required to satisfy
\[
\alpha(e)=UeU^{*},\quad\text{respectively }\Phi(e)=UeU^{*}\quad(e\in E),\qquad\|U-1\|<\kappa,
\]
for some unitary $U\in M$. The bound on $U-1$ is in operator norm.
\end{prop}

\begin{proof}
Write $G=\operatorname{Aut}_{\tau,Q}(M)$ for the group of trace-preserving
automorphisms fixing $Q$ pointwise. We prove (1)$\Rightarrow$(2),
including the additional assertion about $E$; taking $E=\mathbb{C}1$
gives the original assertion. We use the Baire category argument of
\cite[Proposition 1]{ShlyakhtenkoFuchsian}, retaining both the restriction
of an automorphism and its implementing unitary.

Choose $L^{2}$-dense sequences $(a_{j})_{j\geq1}$ in $N$ and $(z_{j})_{j\geq1}$
in $M$. On the space of normal trace-preserving embeddings $N\to M$
fixing $Q$ pointwise, put 
\[
d(\phi,\psi)=\sum_{j=1}^{\infty}2^{-j}\min\{1,\|\phi(a_{j})-\psi(a_{j})\|_{2}\}.
\]
This metric induces pointwise $L^{2}$ convergence. The space is complete.
Indeed, if $(\phi_{n})$ is Cauchy, then for each $a\in N$ the sequence
$(\phi_{n}(a))$ is $L^{2}$-Cauchy and satisfies $\|\phi_{n}(a)\|=\|a\|$.
A uniformly operator-norm bounded $L^{2}$-Cauchy sequence has its
limit in $M$. These limits define a unital linear map on $N$, preserve
products and adjoints by continuity on norm-bounded sets, and preserve
the trace. Faithfulness of the trace makes it injective. It is normal:
if $b_{i}\uparrow b$ in $N_{+}$ and $c=\sup_{i}\phi(b_{i})$, then
$c\leq\phi(b)$ and normality of the traces gives $\tau(c)=\tau(b)=\tau(\phi(b))$,
hence $c=\phi(b)$. The condition of fixing $Q$ is closed.

Let $\mathcal{X}$ be the closure of
\[
\{(\alpha|_{N},U):\alpha\in G,\ U\in\mathcal{U}(M),\ \alpha|_{E}=\operatorname{Ad}U|_{E}\}
\]
in the product of this embedding space and the unitary group with
its operator-norm metric. This is a complete metric space containing
$(\iota,1)$, where $\iota:N\hookrightarrow M$ is inclusion. Every
$(\phi,U)\in\mathcal{X}$ still satisfies $\phi|_{E}=\operatorname{Ad}U|_{E}$.
For $j,k\geq1$, set
\[
\mathcal{O}_{j,k}=\{(\phi,U)\in\mathcal{X}:\text{ there is }y\in N\text{ with }\|z_{j}-\phi(y)\|_{2}<1/k\}.
\]
Each $\mathcal{O}_{j,k}$ is open, being a union of open conditions
on the value at a fixed $y$.

To see density, it suffices to approximate a pair $(\alpha_{0}|_{N},U_{0})$
from the set whose closure defines $\mathcal{X}$. A neighborhood
controls the first coordinate on a finite set $F\subset N$ and the
second coordinate in norm. Apply the strengthened condition (1) to
$F$ and $\alpha_{0}^{-1}(z_{j})$, obtaining $\beta\in G$, $y\in N$,
and an implementing unitary $V$ with both the displacement on $F$
and $\|V-1\|$ arbitrarily small, and with approximation error less
than $1/k$. Then $\alpha_{0}\beta$ is implemented on $E$ by $\alpha_{0}(V)U_{0}$,
and
\[
\begin{aligned}\|\alpha_{0}\beta(a)-\alpha_{0}(a)\|_{2} & =\|\beta(a)-a\|_{2}\quad(a\in F),\\
\|\alpha_{0}(V)U_{0}-U_{0}\| & =\|V-1\|,\\
\|z_{j}-\alpha_{0}\beta(y)\|_{2} & <1/k.
\end{aligned}
\]
Thus the resulting pair belongs to the given neighborhood and to $\mathcal{O}_{j,k}$.

By the Baire category theorem, $\bigcap_{j,k}\mathcal{O}_{j,k}$ is
dense in $\mathcal{X}$. If $(\phi,U)$ belongs to this intersection,
each $z_{j}$ belongs to the $L^{2}$ closure of $\phi(N)$. The image
$\phi(N)$ is a von Neumann subalgebra of $M$. Its trace-preserving
conditional expectation therefore fixes every $z_{j}$. By $L^{2}$
density, this expectation is the identity, so $\phi(N)=M$. Choosing
$(\phi,U)$ sufficiently close to $(\iota,1)$ on the prescribed finite
set and in the unitary coordinate proves (2), with the additional
bound on $\|U-1\|$.

For the original conditions, assume (2). By separability, choose surjective
trace-preserving isomorphisms $\Phi_{j}:N\to M$, fixing $Q$, which
converge pointwise in $L^{2}$ to the inclusion. Fix one surjective
trace-preserving isomorphism $\Phi:N\to M$ fixing $Q$. With $\iota:N\hookrightarrow M$
denoting inclusion, the maps 
\[
\gamma_{j}=\Phi_{j}\Phi^{-1}\in G
\]
converge pointwise in $L^{2}$ to $\rho=\iota\Phi^{-1}$, whose range
is $N$. This proves (3). Conversely, if $\gamma_{j}\to\rho$ as in
(3), then 
\[
\Phi_{j}=\gamma_{j}\rho^{-1}:N\to M
\]
are surjective trace-preserving isomorphisms fixing $Q$ and converging
to the inclusion. Thus (3)$\Rightarrow$(2).

It remains to prove (2)$\Rightarrow$(1), again including the additional
assertion. Given $F,z,\eta,\epsilon,\kappa$, choose a surjective
isomorphism $\Phi:N\to M$ fixing $Q$, moving each $a\in F$ by less
than $\eta/2$, and implemented on $E$ by a unitary $U$ with $\|U-1\|<\kappa/2$.
Put $y=\Phi^{-1}(z)\in N$. Choose a second surjective isomorphism
$\Psi:N\to M$ fixing $Q$, with
\[
\max_{a\in F}\|\Psi(a)-a\|_{2}<\eta/2,\qquad\|\Psi(y)-y\|_{2}<\epsilon,
\]
and implemented on $E$ by a unitary $V$ with $\|V-1\|<\kappa/2$.
The automorphism $\alpha=\Phi\Psi^{-1}$ belongs to $G$, and
\[
\|\alpha(a)-a\|_{2}\leq\|a-\Psi(a)\|_{2}+\|\Phi(a)-a\|_{2}<\eta\quad(a\in F),
\]
while
\[
\|\alpha(y)-z\|_{2}=\|\Psi^{-1}(y)-y\|_{2}=\|y-\Psi(y)\|_{2}<\epsilon.
\]
For $e\in E$, the identity $\alpha(VeV^{*})=UeU^{*}$ shows that
$\alpha|_{E}$ is implemented by $\alpha(V^{*})U$, and
\[
\|\alpha(V^{*})U-1\|\leq\|V-1\|+\|U-1\|<\kappa.
\]
This proves (1).
\end{proof}
We state a $C^{*}$-algebra version of Proposition \ref{prop:aap}.
The same argument applies in operator norm: point-norm limits of embeddings
remain isometric, and their ranges are norm closed.
\begin{prop}
\label{prop:caap}Let $A\subset B\subset C$ be unital inclusions
of separable $C^{*}$-algebras, and let $\tau$ be a trace on $C$.
The following conditions are equivalent:
\begin{enumerate}
\item For every finite subset $F\subset B$, every $z\in C$, and every
$\eta,\epsilon>0$, there are a trace-preserving automorphism $\alpha$
of $C$ fixing $A$ pointwise and $y\in B$ such that 
\[
\max_{b\in F}\|\alpha(b)-b\|<\eta,\qquad\|z-\alpha(y)\|<\epsilon.
\]
\item For every finite subset $F\subset B$ and every $\epsilon>0$, there
is a surjective trace-preserving isomorphism $\Phi:B\to C$ fixing
$A$ pointwise such that 
\[
\max_{b\in F}\|\Phi(b)-b\|<\epsilon.
\]
\item There is a trace-preserving endomorphism $\rho:C\to C$ with $\rho(C)=B$
which is a point-norm limit of a sequence of trace-preserving automorphisms
of $C$ fixing $A$ pointwise.
\end{enumerate}
Moreover, for a specified unital $C^{*}$-subalgebra $E\subset B$,
conditions (1) and (2) remain equivalent with the additional requirement
that, for every $\kappa>0$, the map be implemented on $E$ by a unitary
$U\in C$ with $\|U-1\|<\kappa$.
\end{prop}

One can consider a version of CAAP where no requirement of preserving
traces (or any particular states) is made; we will not explore that
definition in the present paper.
\begin{defn}
\label{def:aap}An inclusion $N\subset M$ as in Proposition \ref{prop:aap}
has the Automorphism Approximation Property (AAP) relative to $Q\subset N$
if the equivalent conditions of that proposition hold. When $Q=\mathbb{C}1$,
we simply say that $N\subset M$ has AAP.
\end{defn}

\begin{defn}
\label{def:caap}An inclusion $B\subset C$ as in Proposition \ref{prop:caap}
has the $C^{*}$ Automorphism Approximation Property (CAAP) relative
to $A\subset B$ if the equivalent conditions of that proposition
hold. When $A=\mathbb{C}1$, we simply say that $B\subset C$ has
CAAP.
\end{defn}

\begin{defn}
\label{def:inner-implementation}We say that AAP, respectively CAAP,
\textit{holds with small unitary implementation on a specified unital
subalgebra $E$} of the smaller algebra if the strengthened conditions
(1) and (2) in Proposition \ref{prop:aap}, respectively Proposition
\ref{prop:caap}, hold. The subalgebra $E$ is a von Neumann algebra
in the AAP case and a $C^{*}$-algebra in the CAAP case. Taking $E=\mathbb{C}1$
recovers the ordinary property.
\end{defn}

\begin{rem}
\label{rem:closures}(1) Suppose that $A\subset B\subset C$ are separable
unital $C^{*}$-algebras, that $C$ carries a faithful tracial state
$\tau$, and that $B\subset C$ has CAAP relative to $A$, with trace-preserving
witnessing maps. Then, in the GNS representation of $\tau$, the inclusion
$W^{*}(B)\subset W^{*}(C)$ has AAP relative to $W^{*}(A)$. If CAAP
holds with small unitary implementation on $E\subset B$, then AAP
holds with small unitary implementation on $W^{*}(E)$: extend the
surjective trace-preserving isomorphisms to the generated von Neumann
algebras, retaining the same implementing unitaries.

(2) Having AAP or CAAP relative to a subalgebra implies the same property
relative to any smaller unital subalgebra.
\end{rem}

\begin{rem}
\label{rem:implementation}If $N=W^{*}(Q,E)$, then fixing $Q$ and
conjugating $E$ by unitaries converging in norm to $1$ already imply
pointwise $L^{2}$ convergence to the inclusion. The corresponding
statement holds in norm when $B=C^{*}(A,E)$. Thus, in this situation,
the finite-set approximation condition follows from the requirement
of small unitary implementation.

For CAAP, small unitary implementation on $E$ is equivalent to small
unitary implementation on $TET^{*}$ for any $T\in\mathcal{U}(B)$,
with the fixed subalgebra $A$ unchanged. Indeed, if $\Phi|_{E}=\operatorname{Ad}U|_{E}$,
then $\Phi$ is implemented on $TET^{*}$ by $W=\Phi(T)UT^{*}$, and
\[
\|W-1\|\leq\|\Phi(T)-T\|+\|U-1\|.
\]
We can choose $F$ to contain $T$, thus bounding $\Vert\Phi(T)-T\Vert$.
Replacing $T$ by $T^{*}$ proves the converse.
\end{rem}

\subsection{Examples of inclusions with AAP and CAAP}

The simplest examples of inclusions satisfying AAP come from infinite
free products or tensor products:
\begin{example}
\label{ex:infinite-free-product}Let $(P_{j},\tau_{j})$, $j\geq0$,
be trace-preserving copies of a tracial von Neumann algebra $(P,\tau)$
with separable predual. Put $N=*_{j\geq1}(P_{j},\tau_{j})$ and $M=(P_{0},\tau_{0})*N$.
Then $N\subset M$ has AAP, relative to any fixed initial free product
$Q=P_{1}*\cdots*P_{r}$.

Similarly, let $(\mathcal{B},\tau)$ be a separable unital $C^{*}$-algebra
with a faithful tracial state, and let $(B_{j},\tau_{j})$, $j\geq0$,
be trace-preserving copies of $(\mathcal{B},\tau)$. Then the reduced
free-product inclusion $B=*_{\mathrm{red},\,j\geq1}(B_{j},\tau_{j})\subset(B_{0},\tau_{0})*_{\mathrm{red}}B$
has CAAP relative to $A=*_{\mathrm{red},\,j=1}^{r}(B_{j},\tau_{j})$
for every $r\geq1$.
\end{example}

\begin{proof}
Given a finite $F\subset N$ and $z\in M$ and $\eta,\epsilon>0$,
choose $m\geq r$ so that every $a\in F$ has an approximation $a^{\prime}\in P_{1}*\cdots*P_{m}$
with $\|a-a^{\prime}\|_{2}<\eta/2$, and choose $z^{\prime}\in P_{0}*P_{1}*\cdots*P_{m}$
with $\|z-z^{\prime}\|_{2}<\epsilon$. Choose $j>m$. Consider the
automorphism $\alpha:M\to M$ that switches $P_{0}$ and $P_{j}$,
and fixes all other components. Then $\alpha$ fixes $P_{1}*\cdots*P_{m}\supset Q$
pointwise. Moreover 
\[
P_{0}*P_{1}*\cdots*P_{m}\subset\alpha(N).
\]
Hence $y=\alpha^{-1}(z^{\prime})$ belongs to $N$, and 
\[
\|\alpha(a)-a\|_{2}\leq2\|a-a^{\prime}\|_{2}<\eta\quad(a\in F),\qquad\|z-\alpha(y)\|_{2}<\epsilon.
\]
This proves AAP relative to $Q$. The proof of the $C^{*}$-algebra
version is identical.
\end{proof}
\begin{example}
(a) Let $R=\overline{\bigotimes}M_{2\times2}$ be the hyperfinite
II$_{1}$ factor. Then $R\otimes1\subset R\otimes R$ has AAP.

\label{ex:abelian-coordinate}(b) Let $X=\{0,1\}^{\mathbb{N}}$ be
the Cantor set with the product measure $\prod(\frac{1}{2}\delta_{0}+\frac{1}{2}\delta_{1})$.
Then the inclusion 
\[
C(X)\otimes1\subset C(X\times X)
\]
has trace-preserving CAAP.

(c) The coordinate inclusion 
\[
N=L^{\infty}([0,1])\otimes1\subset M=L^{\infty}([0,1]^{2}),
\]
with Lebesgue probability measure, has AAP.
\end{example}

The main results of this paper concern the following examples:
\begin{example}
\label{exa:MainExamples}(a) Let $n\geq2$, $B=C_{\mathrm{red}}^{*}(\mathbb{F}_{n})$
and $C=B*_{\mathrm{red}}C^{*}(S)$, where $S$ is a free semicircular
element. The inclusion $B\subset C$ has CAAP with trace-preserving
witnessing maps. In particular, $B\cong C$. This example is the direct
application of Theorem \ref{thm:free-absorption} proved in this paper.

(b) Let $n\geq2$, $N=L(\mathbb{F}_{n})$ and $M=N*L(\mathbb{Z})\cong L(\mathbb{F}_{n+1})$.
The inclusion $N\subset M$ has AAP, and hence $N\cong M$. This follows
from (a) and Remark \ref{rem:closures}(1), since the von Neumann
algebra generated by a semicircular element is diffuse abelian. 

(c) Let $M=L(\mathbb{F}_{2})$ be generated by free Haar unitaries
$a,b$. Put $w=aba^{-1}b^{-1}$, $N=W^{*}(a,w)$ and $Q=W^{*}(w)$.
The construction in \cite{ShlyakhtenkoFuchsian} gives AAP for $N\subset M$
relative to $Q$. Since $a$ and $w$ are free Haar unitaries, a surjective
isomorphism $\Phi:N\to M$ fixing $Q$ gives 
\[
M=\Phi(W^{*}(w)*W^{*}(a))=W^{*}(w)*W^{*}(v),\qquad v=\Phi(a).
\]
Thus $Q$ is freely complemented in $M$.
\end{example}

\subsection{Further properties of AAP/CAAP inclusions}

The following permanence properties will also be useful later. All
relative properties are with respect to the same fixed unital subalgebra
$Q$.
\begin{prop}
\label{prop:permanence}The following statements hold for AAP in tracial
von Neumann algebras with separable predual, and for CAAP in separable
unital $C^{*}$-algebras with a fixed trace.
\begin{enumerate}
\item If $Q\subset L\subset N\subset M$ and both $L\subset N$ and $N\subset M$
have the property relative to $Q$, then $L\subset M$ has the property
relative to $Q$.
\item Let $Q\subset N_{0}\subset N_{1}\subset\cdots$ be an increasing sequence,
with each $N_{n}\subset N_{n+1}$ having the property relative to
$Q$. Put 
\[
M=\overline{\bigcup_{n\geq0}N_{n}}^{\|\cdot\|}\quad\text{in the }C^{*}\text{ case},\qquad M=\overline{\bigcup_{n\geq0}N_{n}}^{WOT}\quad\text{in the tracial case}.
\]
Then every $N_{n}\subset M$ has the same property relative to $Q$.
\end{enumerate}
Both assertions remain valid with small unitary implementation on
a common unital subalgebra $E$ contained in the smallest algebra.
\end{prop}

\begin{proof}
Use the characterization by surjective isomorphisms in Propositions
\ref{prop:aap}(2) and \ref{prop:caap}(2), including its strengthened
version. Write $d(x,y)=\|x-y\|$ in the $C^{*}$ case and $d(x,y)=\|x-y\|_{2}$
in the tracial case. Taking $E=\mathbb{C}1$ gives the ordinary assertions.

For (1), given a finite $F\subset L$ and $\epsilon,\kappa>0$, choose
a surjective isomorphism $\phi:L\to N$ moving $F$ by less than $\epsilon/2$,
implemented on $E$ by $U$ with $\|U-1\|<\kappa/2$. Choose a surjective
isomorphism $\psi:N\to M$ moving $\phi(F)$ by less than $\epsilon/2$,
implemented on $E$ by $V$ with $\|V-1\|<\kappa/2$. Both maps fix
$Q$. Their composition has the required finite-set approximation
and is implemented on $E$ by $\psi(U)V$, with
\[
\|\psi(U)V-1\|\leq\|U-1\|+\|V-1\|<\kappa.
\]

For (2), it suffices to treat $N_{0}\subset M$. Fix a finite $F\subset N_{0}$
and $\epsilon,\kappa>0$, and choose positive numbers $\epsilon_{n},\kappa_{n}$
with $\sum_{n\geq0}\epsilon_{n}<\epsilon$ and $\sum_{n\geq0}\kappa_{n}<\kappa$.
Choose increasing finite sets $H_{n}\subset N_{n}$ whose union is
$d$-dense in $M$, and a $d$-dense sequence $(a_{j})$ in $N_{0}$.
Starting with $\phi_{0}=\mathrm{id}_{N_{0}}$ and $U_{0}=1$, construct
surjective isomorphisms $\phi_{n}:N_{0}\to N_{n}$ fixing $Q$ and
satisfying $\phi_{n}|_{E}=\operatorname{Ad}U_{n}|_{E}$. At stage
$n$, choose a finite set $F_{n}\subset N_{0}$ containing $F$, all
previous $F_{j}$, the first $n$ elements $a_{j}$, and $\phi_{n}^{-1}(H_{n})$.
Choose a surjective isomorphism $\psi_{n}:N_{n}\to N_{n+1}$ fixing
$Q$, moving every element of $\phi_{n}(F_{n})\cup\{U_{n}\}$ by less
than $\epsilon_{n}$, and implemented on $E$ by $V_{n}$ with $\|V_{n}-1\|<\kappa_{n}$.
Set
\[
\phi_{n+1}=\psi_{n}\phi_{n},\qquad U_{n+1}=\psi_{n}(U_{n})V_{n}.
\]
Then $\phi_{n+1}|_{E}=\operatorname{Ad}U_{n+1}|_{E}$, and
\[
\begin{aligned}d(\phi_{n+1}(a),\phi_{n}(a)) & <\epsilon_{n}\quad(a\in F_{n}),\\
d(U_{n+1},U_{n}) & <\epsilon_{n}+\kappa_{n},\\
\|U_{n+1}-1\| & \leq\|U_{n}-1\|+\kappa_{n}.
\end{aligned}
\]

The maps $\phi_{n}$ converge pointwise for $d$ to an embedding $\phi:N_{0}\to M$,
fixing $Q$ and moving $F$ by less than $\epsilon$. In the tracial
case, the uniform bounds $\|\phi_{n}(a)\|=\|a\|$ ensure that the
limits lie in $M$ and preserve multiplication; trace preservation
gives normality, as in Proposition \ref{prop:aap}. The unitaries
$U_{n}$ converge for $d$ to a unitary $U$. In the tracial case
this follows because both $U_{n}$ and $U_{n}^{*}$ are $L^{2}$-Cauchy
and uniformly bounded. The uniform operator-norm bound passes to the
limit, giving
\[
\phi|_{E}=\operatorname{Ad}U|_{E},\qquad\|U-1\|\leq\sum_{n\geq0}\kappa_{n}<\kappa.
\]

To see surjectivity, fix $b\in\bigcup_{n}H_{n}$. For every sufficiently
large $n$, put $a_{n,b}=\phi_{n}^{-1}(b)\in F_{n}$. Retaining this
element in all subsequent control sets gives
\[
d(\phi(a_{n,b}),b)\leq\sum_{j\geq n}\epsilon_{j}\longrightarrow0.
\]
Thus the range is $d$-dense. In the $C^{*}$ case it is norm closed.
In the tracial case its trace-preserving conditional expectation fixes
a $d$-dense set and hence is the identity. Therefore $\phi$ is surjective.
All maps preserve the specified trace.
\end{proof}

\subsection{Outline of the proof}

The key $C^{*}$-algebra result, Theorem \ref{thm:free-absorption},
establishes CAAP with small unitary implementation for the inclusion
$C^{*}(u,v)\subset C^{*}(u,v,S)$ where $u,v$ are two free Haar unitaries,
and $S$ is a free semicircular element free from $u,v$. In particular,
CAAP requires that there be an automorphism $\alpha$ so that $\alpha(u)\approx u$,
$\alpha(v)\approx v$ and $\alpha(x)\approx S$ for some $x\in C^{*}(u,v)$. 

The automorphism $\alpha$ is constructed as $\alpha_{\pi/2}$, where
$\alpha_{t}$ is a certain one-parameter group of automorphisms. This
one-parameter group is obtained by the free semicircular construction
of \cite{ShlyakhtenkoMicrostates} (based on Voiculescu's method of
exponentiating ``divergence-free'' vector fields \cite{VoiculescuCyclomorphy}).
This construction is described in \S\ref{sec:ExponentiatingDervations}
and takes as input a densely defined derivation $\delta:C^{*}(u,v)\to L^{2}(C^{*}(u,v)\otimes C^{*}(u,v)^{op})$.
Proposition \ref{prop:absorption} gives the construction of $\alpha$
once a derivation with the required properties is obtained. This is
the same strategy as OpenAI's paper \cite{OpenAIFreeGroupFactors},
although implemented differently; one important technical point in
our argument is that the automorphism fixes $v$ and conjugates $u$
by a unitary that is close to $1$. 

The key requirements on the derivation are: $\delta(u)\approx0$,
$\delta(v)=0$ (these correspond to $\alpha(u)\approx u$ and $\alpha(v)=v$).
Additionally, writing $\xi=\delta^{*}(1\otimes1)$, we require $\delta(\xi)\approx1\otimes1$
(this corresponds to $\alpha(\xi)\approx S$). 

The construction of the derivation is carried out in \S\ref{sec:constructionOfDerivation}
using the same cocycle idea as in \cite{OpenAIFreeGroupFactors}.
Finally, \S\ref{sec:Semicircular-absorption} and \S\ref{sec:free-product-consequences}
are devoted to operator-algebraic consequences of Theorem \ref{thm:free-absorption}.
We end by connecting AAP/CAAP with continuous model theory in \S\ref{sec:elementary-inclusions}.

We explain a classical version of our construction in the Appendix. 

\subsection*{Acknowledgements}

This work is based on \cite{OpenAIFreeGroupFactors} and proofs that
are AI-generated. ChatGPT 6.0 Pro was used in the analysis of OpenAI's
proof, and in preparation of this manuscript. 

The author is grateful to Sorin Popa for several useful discussions.

This research was sponsored in part by the Army Research Office and
was accomplished under Grant Number W911NF-25-1-0075. The views and
conclusions contained in this document are those of the authors and
should not be interpreted as representing the official policies, either
expressed or implied, of the Army Research Office or the U.S. Government.
The U.S. Government is authorized to reproduce and distribute reprints
for Government purposes notwithstanding any copyright notation herein.
Research was also supported in part by NSF grant DMS-2348633.

\section{Exponentiating derivations to flows}\label{sec:ExponentiatingDervations}

Let $(B,\tau)$ be a unital $C^{*}$-algebra with a faithful tracial
state, and let $B_{0}\subset B$ be a norm-dense unital $*$-subalgebra.
We write $B_{0}\otimes B_{0}^{\mathrm{op}}$ for the algebraic tensor
product with the multiplication 
\[
(a\otimes b^{\mathrm{op}})(c\otimes d^{\mathrm{op}})=ac\otimes(db)^{\mathrm{op}}.
\]
Its bimodule structure is given by 
\[
a\cdot(u\otimes v^{\mathrm{op}})\cdot b=au\otimes(vb)^{\mathrm{op}}.
\]
Let $\delta:B_{0}\to B_{0}\otimes B_{0}^{\mathrm{op}}$ be a derivation.
We say that $\delta$ is real if 
\[
\delta(a^{*})=\delta(a)^{\dagger},\qquad(u\otimes v^{\mathrm{op}})^{\dagger}=v^{*}\otimes(u^{*})^{\mathrm{op}}.
\]
Here $\dagger$ is the flip-adjoint
\[
(u\otimes v^{\mathrm{op}})^{\dagger}=v^{*}\otimes(u^{*})^{\mathrm{op}}.
\]

Regard $\delta$ as a densely defined operator from $L^{2}(B,\tau)$
to $L^{2}(B,\tau)\otimes L^{2}(B^{\mathrm{op}},\tau)$, with inner
products linear in the second variable. Suppose that $\delta$ is
real and that $\xi=\delta^{*}(1\otimes1)$ exists and belongs to $B_{0}$.
In that case, $\xi$ is called Voiculescu's free conjugate variable
associated with $\delta$ (see \cite{VoiculescuEntropyV}; for a formulation
involving more general derivations, see the partial conjugate variables
of \cite[Section 2.1, Definition 2.1]{ShlyakhtenkoPartialConjugates}).
Since $\delta$ is real, it follows that $\xi=\xi^{*}$; the equality
$\xi=\delta^{*}(1\otimes1)$ is then 
\[
\tau(\xi b)=(\tau\otimes\tau)(\delta(b)),\qquad b\in B_{0}.
\]
In particular $\tau(\xi)=0$.
\begin{prop}
\label{prop:adjoint}Under these assumptions, $B_{0}\otimes B_{0}^{\mathrm{op}}\subset\operatorname{Dom}(\delta^{*})$,
and for $a,b\in B_{0}$ one has 
\[
\delta^{*}(a\otimes b^{\mathrm{op}})=a\xi b-(\mathrm{id}\otimes\tau)(\delta(a))b-a(\tau\otimes\mathrm{id})(\delta(b)).
\]
Consequently $\delta$ is closable.
\end{prop}

\begin{proof}
This is the standard adjoint calculation of \cite[Corollary 4.2 and Proposition 4.3]{VoiculescuEntropyV};
see also \cite[proof of Lemma 6]{ShlyakhtenkoMicrostates}. Since
$B_{0}\otimes B_{0}^{\mathrm{op}}$ is dense, $\delta$ is closable.
\end{proof}
For an algebraic tensor $U=\sum_{i}a_{i}\otimes b_{i}^{\mathrm{op}}$
and an element $z$ of an algebra containing $B_{0}$, write 
\[
U\#z=\sum_{i}a_{i}zb_{i}.
\]
For tensors $U,T$, we also write $U\#T=UT$, using the tensor-algebra
multiplication above; thus $(U\#T)\#z=U\#(T\#z)$.

Let $S$ be a standard semicircular element, free from $B$; thus
$\Vert S\Vert=2$ and $\Vert S\Vert_{2}=1$.

Put 
\[
C=(B,\tau)*(C^{*}(S),\sigma)=(B,\tau)*_{\mathrm{red}}(C([-2,2]),\sigma),\qquad C_{0}=\mathop{{*}\text{-}\mathrm{Alg}}\nolimits(B_{0},S),
\]
where $\sigma$ is the semicircular probability measure of variance
one. We also write $\tau$ for the faithful free-product trace.
\begin{lem}
\label{lem:flow}Suppose that $B_{0}=\mathop{{*}\text{-}\mathrm{Alg}}\nolimits(A,X)$,
where $A\subset B_{0}$ is a unital $*$-subalgebra and $X\in B_{0}$.
Suppose that $\delta$ is real, $\delta|_{A}=0$, and $\xi=\delta^{*}(1\otimes1)\in B_{0}$.
Then the formulas 
\[
D(b)=\delta(b)\#S\quad(b\in B_{0}),\qquad D(S)=-\xi
\]
define a $*$-derivation $D:C_{0}\to C_{0}$ satisfying $\tau(Dp)=0$
for every $p\in C_{0}$. It generates a global point-norm continuous
one-parameter group $(\alpha_{t})_{t\in\mathbb{R}}$ of trace-preserving
automorphisms of $C$, fixing $A$ pointwise, with 
\[
\frac{d}{dt}\alpha_{t}(p)=\alpha_{t}(Dp),\qquad p\in C_{0}.
\]
In particular, 
\[
\|\alpha_{t}(X)-X\|\leq|t|\,\|\delta(X)\#S\|.
\]

If the restriction of $\delta$ to a unital $*$-subalgebra $E\subset B_{0}$
is inner, say 
\[
\delta(b)=Qb-bQ\quad(b\in E),\qquad Q\in B_{0}\otimes B_{0}^{\mathrm{op}},\qquad Q^{\dagger}=-Q,
\]
then there are norm continuously differentiable unitaries $U_{t}\in C$,
with $U_{0}=1$, such that 
\[
\alpha_{t}(b)=U_{t}bU_{t}^{*}\quad(b\in C^{*}(E)),\qquad\|U_{t}-1\|\leq|t|\,\|Q\#S\|.
\]
\end{lem}

\begin{proof}
The construction in \cite[Appendix, Proposition 21]{ShlyakhtenkoMicrostates}
gives the derivation $D$ and the identity $\tau(Dp)=0$. Only finitely
many coefficients from $A$ occur in the polynomial expressions for
$D(X)$ and $D(S)$; note that $D$ annihilates $A$. The analytic-vector
argument in \cite[Theorem 2.7 and Corollary 3.7]{VoiculescuCyclomorphy}
therefore yields a one-parameter group $(\alpha_{t})$ of trace-preserving
automorphisms of $W^{*}(C)$. For sufficiently small $|t|$, the power
series for $\alpha_{t}(X)$ and $\alpha_{t}(S)$ converge in operator
norm to elements of $C$. Since $\alpha_{t}$ fixes $A$, this gives
$\alpha_{t}(C)\subset C$. Applying the same argument to $-t$ gives
equality for sufficiently small $|t|$. For arbitrary $t\in\mathbb{R}$,
choose $n$ so large that this applies to $t/n$; then $\alpha_{t}=(\alpha_{t/n})^{n}$
gives $\alpha_{t}(C)=C$. The maps $t\mapsto\alpha_{t}(y)$ are norm
continuous for every $y\in C$ and norm differentiable for $y\in C_{0}$.

For the operator-norm estimate needed here, norm differentiability
gives 
\[
\alpha_{t}(y)-y=\int_{0}^{t}\alpha_{s}(D(y))\,ds\qquad\forall y\in C_{0},
\]
with the integral taken in operator norm. Since each $\alpha_{s}$
is isometric, 
\[
\|\alpha_{t}(y)-y\|\leq|t|\,\|D(y)\|.
\]
Taking $y=X$ gives the first estimate.

For the restriction to $E$, put $H=Q\#S$. Then $H^{*}=-H$ and $D(b)=[H,b]$
for $b\in E$. Consider 
\[
\frac{d}{dt}U_{t}=\alpha_{t}(H)U_{t},\qquad U_{0}=1.
\]
This bounded linear differential equation has a global unitary solution.
Both $\alpha_{t}(b)$ and $U_{t}bU_{t}^{*}$ solve $\frac{d}{dt}Y_{t}=[\alpha_{t}(H),Y_{t}]$
with initial value $b$, so uniqueness gives the claimed equality,
first on $E$ and then on $C^{*}(E)$. Finally, integration gives
\[
\|U_{t}-1\|\leq|t|\,\|H\|.
\]
\end{proof}
For $U=\sum_{i}a_{i}\otimes b_{i}^{\mathrm{op}}\in B_{0}\otimes B_{0}^{\mathrm{op}}$,
define the norm 
\[
\|U\|_{\mathrm{rc}}=\bigl\|\sum_{i,j}b_{i}^{*}\tau(a_{i}^{*}a_{j})b_{j}\bigr\|^{1/2}+\bigl\|\sum_{i,j}a_{i}\tau(b_{i}b_{j}^{*})a_{j}^{*}\bigr\|^{1/2}.
\]
The following lemma relates this norm to $\|U\#S\|$.
\begin{lem}
\label{lem:gaussian-bound}Let $S$ be standard semicircular and free
from a faithfully tracial $C^{*}$-algebra $B$. For a finite tensor
$U\in B\otimes B^{\mathrm{op}}$, 
\[
\frac{1}{2}\|U\|_{\mathrm{rc}}\leq\|U\#S\|\leq\|U\|_{\mathrm{rc}}.
\]
\end{lem}

\begin{proof}
Write $U=\sum_{i}a_{i}\otimes b_{i}^{\mathrm{op}}$. Use the Fock
realization over the coarse bimodule of $M=W^{*}(B)$ described in
Connes-{}-Shlyakhtenko \cite[Appendix II]{ConnesShlyakhtenko}. In
this realization, $S=\ell+\ell^{*}$ and $\ell^{*}b\ell=\tau(b)1$
for $b\in M$; see \cite[equation (6.1)]{ConnesShlyakhtenko}. For
$L=\sum_{i}a_{i}\ell b_{i}$, one has 
\[
L^{*}L=\sum_{i,j}b_{i}^{*}\tau(a_{i}^{*}a_{j})b_{j},
\]
so that $\|L\|=\bigl\|\sum_{i,j}b_{i}^{*}\tau(a_{i}^{*}a_{j})b_{j}\bigr\|^{1/2}$.
Applying the same calculation to the adjoint of $\sum_{i}a_{i}\ell^{*}b_{i}$
and using the triangle inequality proves the upper bound.

For the lower bound, put $Y=U\#S$ and let $E_{M}$ be the trace-preserving
conditional expectation onto $M$. Freeness gives 
\[
E_{M}(Y^{*}Y)=\sum_{i,j}b_{i}^{*}\tau(a_{i}^{*}a_{j})b_{j},\qquad E_{M}(YY^{*})=\sum_{i,j}a_{i}\tau(b_{i}b_{j}^{*})a_{j}^{*}.
\]
Since $E_{M}$ is contractive, each term in $\|U\|_{\mathrm{rc}}$
is at most $\|Y\|$.
\end{proof}
\begin{prop}
\label{prop:absorption}Let $(B,\tau)$ be separable, unital and faithfully
tracial. Let $B_{0}\subset B$ be a dense unital $*$-subalgebra algebraically
generated by a unital $*$-subalgebra $A$ and an element $X$. Let
$S$ and $C$ be as above. Assume that for every $\epsilon,\eta>0$
there is a real derivation $\delta:B_{0}\to B_{0}\otimes B_{0}^{\mathrm{op}}$
such that:
\begin{enumerate}
\item $\delta|_{A}=0$;
\item The conjugate variable $\xi=\delta^{*}(1\otimes1)$ exists and belongs
to $B_{0}$;
\item $\|\delta(X)\|_{\mathrm{rc}}<\eta$;
\item $\|\delta(\xi)-1\otimes1\|_{\mathrm{rc}}<\epsilon$.
\end{enumerate}
Then $B\subset C$ has CAAP relative to $C^{*}(A)$. Thus $N=W^{*}(B)\subset W^{*}(C)=M$
has AAP relative to $Q=W^{*}(A)$.

Moreover, suppose that one of these derivations is inner on a unital
$*$-subalgebra $E\subset B_{0}$, with $\delta(b)=\Theta b-b\Theta$
for $b\in E$, where $\Theta\in B_{0}\otimes B_{0}^{\mathrm{op}}$
and $\Theta^{\dagger}=-\Theta$. Then the local approximating automorphism
$\alpha=\alpha_{\pi/2}$ constructed in the proof satisfies
\[
\alpha(b)=UbU^{*}\quad(b\in C^{*}(E)),\qquad\|U-1\|\leq\frac{\pi}{2}\|\Theta\#S\|\leq\frac{\pi}{2}\|\Theta\|_{\mathrm{rc}}
\]
for a unitary $U\in C$. Consequently, if for every $\epsilon,\eta,\kappa>0$
the derivation can be chosen to satisfy (1)-{}-(4) and to be inner
on this fixed $E$ with $\|\Theta\|_{\mathrm{rc}}<\kappa$, then CAAP
holds with small unitary implementation on $C^{*}(E)$, and AAP holds
with small unitary implementation on $W^{*}(E)$.
\end{prop}

\begin{proof}
Let $\delta,\xi$ satisfy the hypotheses, and let $(\alpha_{t})$
be the one-parameter group supplied by Lemma \ref{lem:flow}. Put
\[
e=D(\xi)-S=(\delta(\xi)-1\otimes1)\#S.
\]
By Lemma \ref{lem:gaussian-bound}, $\|e\|\leq\|\delta(\xi)-1\otimes1\|_{\mathrm{rc}}<\epsilon$,
and 
\[
\frac{d}{dt}\alpha_{t}(\xi)=\alpha_{t}(S)+\alpha_{t}(e),\qquad\frac{d}{dt}\alpha_{t}(S)=-\alpha_{t}(\xi).
\]
Thus
\[
\alpha_{t}(\xi)=(\cos t)\xi+(\sin t)S+\int_{0}^{t}\cos(t-s)\alpha_{s}(e)\,ds.
\]
Since every $\alpha_{s}$ is isometric, evaluating at $t=\pi/2$ gives
\begin{equation}
\|\alpha_{\pi/2}(\xi)-S\|\leq\|e\|<\epsilon.\label{eq:xiMovedToS}
\end{equation}
At the same time, $\|\delta(X)\#S\|\leq\|\delta(X)\|_{\mathrm{rc}}<\eta$,
so 
\[
\|\alpha_{\pi/2}(X)-X\|<\frac{\pi}{2}\eta,\qquad\alpha_{\pi/2}|_{C^{*}(A)}=\mathrm{id}.
\]

The estimate for inner restrictions follows from Lemma \ref{lem:flow},
with $U=U_{\pi/2}$, and Lemma \ref{lem:gaussian-bound}. Under the
additional hypothesis, we can impose $\|\Theta\|_{\mathrm{rc}}<2\kappa/\pi$
along with the choices of $\eta$ and $\epsilon$ below. The resulting
automorphism then has an implementing unitary within $\kappa$ of
$1$, so the same argument verifies the strengthened condition (1)
of Proposition \ref{prop:caap}.

We now verify condition (1) of Proposition \ref{prop:caap}. Let $F=\{z_{1},\dots,z_{k}\}\subset B$
be a finite set, $z\in C$, and $\eta^{\prime},\epsilon^{\prime}>0$.
Choose $*$-polynomials $p_{1},\dots,p_{k}$ with coefficients from
$A$ so that $\Vert p_{j}(X)-z_{j}\Vert<\eta'/3$, $j=1,\dots,k$,
and a $*$-polynomial in two noncommuting variables with coefficients
from $A$ so that $\Vert p(X,S)-z\Vert<\epsilon'/2$. These polynomial
functions are Lipschitz on norm-bounded sets. Taking $\epsilon<1$,
the preceding estimates give $\|\alpha_{\pi/2}(\xi)\|<3$ and $\|\alpha_{\pi/2}(X)\|=\|X\|$,
so the Lipschitz constants may be chosen independently of the derivation.
We can therefore choose $\eta$ and $\epsilon$ sufficiently small
so that, writing $\alpha=\alpha_{\pi/2}$,
\[
\Vert p_{j}(\alpha(X))-p_{j}(X)\Vert<\eta'/3,\qquad j=1,\dots,k
\]
and
\[
\Vert p(\alpha(X),\alpha(\xi))-p(X,S)\Vert<\epsilon'/2.
\]
Noting that $\alpha(p_{j}(X))=p_{j}(\alpha(X))$ gives us
\[
\begin{aligned}\|\alpha(z_{j})-z_{j}\| & \leq\|\alpha(z_{j})-\alpha(p_{j}(X))\|\\
 & \quad+\|p_{j}(\alpha(X))-p_{j}(X)\|+\|p_{j}(X)-z_{j}\|<\eta'.
\end{aligned}
\]
Likewise, noting that $\alpha(p(X,\xi))=p(\alpha(X),\alpha(\xi))$,
and setting $y=p(X,\xi)\in B$,
\[
\|z-\alpha(y)\|\leq\|z-p(X,S)\|+\|p(X,S)-p(\alpha(X),\alpha(\xi))\|<\epsilon'.
\]
Thus condition (1) is verified, and the proof of CAAP is finished.
Remark \ref{rem:closures}(1) now gives AAP for $N\subset M$ relative
to $Q$.
\end{proof}
\begin{rem}
\label{rem:xiIsCloseToSemicircular}Equation \ref{eq:xiMovedToS},
with $\alpha=\alpha_{\pi/2}$, gives 
\[
\|\xi-\alpha^{-1}(S)\|=\|\alpha(\xi)-S\|<\epsilon.
\]
Since $\alpha$ preserves the trace, $\alpha^{-1}(S)$ is standard
semicircular. Thus $\xi$ is within $\epsilon$ in operator norm of
a semicircular element in the enlarged algebra $C$.
\end{rem}

\section{Constructing a liberation derivation}\label{sec:constructionOfDerivation}

Let $G=\mathbb{F}_{2}=\langle a,b\rangle$. The goal of this section
is to construct the derivations required by Proposition \ref{prop:absorption}
for 
\[
B=C_{\mathrm{red}}^{*}(G),\qquad B_{0}=\mathbb{C}G,\qquad A=\mathop{{*}\text{-}\mathrm{Alg}}\nolimits(u_{b}),\qquad X=u_{a},
\]
with the canonical trace on $B$. Thus $B_{0}=\mathop{{*}\text{-}\mathrm{Alg}}\nolimits(A,X)$.
For every $\eta,\epsilon>0$, we seek a real derivation $\delta:B_{0}\to B_{0}\otimes B_{0}^{\mathrm{op}}$
such that 
\[
\delta|_{A}=0,\qquad\xi=\delta^{*}(1\otimes1)\in B_{0},\qquad\|\delta(u_{a})\|_{\mathrm{rc}}<\eta,\qquad\|\delta(\xi)-1\otimes1\|_{\mathrm{rc}}<\epsilon.
\]

As discussed in Remark \ref{rem:xiIsCloseToSemicircular}, $\xi$
must then be close in operator norm to a semicircular element in the
enlarged algebra.

We first use a group cocycle to construct a derivation with small
values at the generators. We then choose a polynomial kernel from
a regularized semicircular change of variables and make a small further
correction so that the full ambient conjugate variable is exactly
the prescribed polynomial $p(u_{w})$.

For any discrete group $\Gamma$, let $(\delta_{g})_{g\in\Gamma}$
be the standard orthonormal basis of $\ell^{2}(\Gamma)$, and identify
$\delta_{g}$ with the canonical group unitary $u_{g}$ in $L^{2}(L(\Gamma))$.
Throughout, $\Delta$ denotes the group comultiplication in the outer-bimodule
convention: 
\[
\begin{aligned}\Delta:\ell^{2}(\Gamma) & \longrightarrow\ell^{2}(\Gamma)\otimes\ell^{2}(\Gamma),\\
\Delta(\delta_{g}) & =\delta_{g}\otimes\delta_{g^{-1}}.
\end{aligned}
\]
If $\lambda$ denotes the left regular representation, then 
\[
\Delta(\lambda(g)f)=u_{g}\Delta(f)u_{g}^{*},\qquad\Delta(\delta_{e})=1\otimes1.
\]

We start by proving Lemma \ref{lem:cocycle}, which supplies an element
$x\in B_{0}\subset\ell^{2}(G)$, a group element $w\in G\setminus\{e\}$
and a cocycle $c:G\to\ell^{2}(G)$ uniquely determined on generators
by 
\[
c(a)=(1-\lambda(a))x,\qquad c(b)=0,
\]
for which both 
\[
\rho=\|x\|_{2},\qquad\theta=\|c(w)-\delta_{e}\|_{2}
\]
are arbitrarily small. The idea is to take $w$ to be a long word
consisting of long runs of $a$ and $b$. Applying the cocycle identity
writes $c(w)$ as a sum of group elements which are essentially freely
independent from each other, applied to $x$. Using the free central
limit theorem, we then argue that it is possible to first choose such
a $w$ and then choose $x$ in such a way that the cocycle $c(a)\approx0$,
$c(b)=0$, but $c(w)\approx\delta_{e}\in\ell^{2}(G)$.

The cocycle identity and the equivariance of $\Delta$ define an algebraic
derivation 
\[
\partial(u_{g})=\Delta(c(g))u_{g}\qquad(g\in G).
\]
It satisfies 
\[
\partial(u_{a})\approx0,\qquad\partial(u_{b})=0,\qquad\partial(u_{w})\approx1\otimes u_{w}^{\mathrm{op}}.
\]
Consider the corrected derivation $\hat{\partial}$ given by
\[
\widehat{\partial}(z)=c(w)(e)^{-1}\partial(z)\#(1\otimes(u_{w}^{*})^{\mathrm{op}})
\]
(note that $c(w)\approx\delta_{e}$ so $c(w)(e)\approx1$). Then $\widehat{\partial}(u_{w})\approx1\otimes1$,
$\hat{\partial}(u_{a})\approx0$, $\hat{\partial}(u_{b})=0$ and 
\[
\widehat{\partial}^{\,*}(1\otimes1)=u_{w}\in\mathop{{*}\text{-}\mathrm{Alg}}\nolimits(u_{w}).
\]
However, $\hat{\partial}$ does not give us the derivation we want,
because it is not real (in particular, $\hat{\partial}^{*}(1\otimes1)$
is not self-adjoint).

To motivate the rest of the construction, we proceed formally, ignoring
domain issues and pretending that instead of approximate equality
we have the exact identity $\hat{\partial}(u_{w})=1\otimes1$.

For the desired derivation $\delta$, $\delta^{*}(1\otimes1)$ should
be close to a semicircular element (Remark \ref{rem:xiIsCloseToSemicircular}).
Thus we choose a continuous function $f$ such that $Y=f(u_{w})$
is semicircular.

Let $\mathscr{A}=W^{*}(u_{w})$. Identify $\mathscr{A}\cong L^{\infty}(\mathbb{T})$
with Haar measure, $\mathscr{A}\bar{\otimes}\mathscr{A}^{\mathrm{op}}$
with $L^{\infty}(\mathbb{T}\times\mathbb{T})$, and denote by $\iota:\mathscr{A}\bar{\otimes}\mathscr{A}^{\mathrm{op}}\to L^{2}(L(G)\bar{\otimes}L(G)^{\mathrm{op}})$
the inclusion map. Under this identification, formally, 
\[
\widehat{\partial}(g(u_{w}))=\iota\!\left(\frac{g(\zeta)-g(\omega)}{\zeta-\omega}\right).
\]
Thus 
\[
\widehat{\partial}(Y)=\iota(Df),\qquad Df=\frac{f(\zeta)-f(\omega)}{\zeta-\omega}.
\]
Let 
\[
\widehat{K}(\zeta,\omega)=\frac{\zeta-\omega}{f(\zeta)-f(\omega)}.
\]
Setting $\delta(z)=\widehat{\partial}(z)\#\iota(\widehat{K})$ gives
\[
\delta(Y)=\iota(Df)\#\iota(\widehat{K})=\iota\bigl((Df)\widehat{K}\bigr)=1\otimes1.
\]
Equivalently, $\delta(z)=\partial(z)\#\iota(K)$, where 
\[
K(\zeta,\omega)=\frac{1}{c(w)(e)}\,\frac{\zeta-\omega}{f(\zeta)-f(\omega)}\,\frac{1}{\omega}.
\]

The crucial point is that $\widehat{\partial}^{\,*}(1\otimes1)=u_{w}\in\mathscr{A}$,
and $\widehat{\partial}$ takes $\mathscr{A}$ into its tensor square.
The adjoint calculation underlying Proposition \ref{prop:adjoint}
gives formally 
\[
\widehat{\partial}^{\,*}\bigl(\iota(\mathscr{A}\otimes\mathscr{A}^{\mathrm{op}})\bigr)\subset L^{2}(\mathscr{A}).
\]
Since the correcting kernel $\widehat{K}$ belongs formally to the
tensor square of the same abelian algebra $\mathscr{A}$, we get 
\[
\delta^{*}(1\otimes1)=\widehat{\partial}^{\,*}\bigl(\iota(\widehat{K}^{*})\bigr)\in L^{2}(\mathscr{A}).
\]
Thus the conjugate variable computed inside $\mathscr{A}$ is already
the full ambient conjugate variable.

Let $\mathscr{A}_{0}=\mathop{{*}\text{-}\mathrm{Alg}}\nolimits(Y)$.
The restriction of $\delta$ to $\mathscr{A}_{0}$ is the composition
of Voiculescu's difference quotient $\mathcal{D}_{Y}:\mathscr{A}_{0}\to\mathscr{A}_{0}\otimes\mathscr{A}_{0}^{\mathrm{op}}$,
defined by $\mathcal{D}_{Y}(Y)=1\otimes1$, with the two inclusions
\[
\mathscr{A}_{0}\otimes\mathscr{A}_{0}^{\mathrm{op}}\subset\mathscr{A}\bar{\otimes}\mathscr{A}^{\mathrm{op}}\xrightarrow{\ \iota\ }L^{2}(L(G)\bar{\otimes}L(G)^{\mathrm{op}}).
\]
The scalar adjoint calculation in $\mathscr{A}$ now gives 
\[
\delta^{*}(1\otimes1)=\left[2\,\mathrm{p.v.}\!\int_{\mathbb{T}}\frac{1}{f(\zeta)-f(\omega)}\,dm(\omega)\right]_{\zeta=u_{w}}=f(u_{w})=Y,
\]
since $f_{*}m$ is the semicircular law \cite{VoiculescuEntropyV}.
Hence $\xi=Y$ and $\delta(\xi)=1\otimes1$, as desired.

The rigorous proof regularizes this kernel, approximates the regularized
kernel by Laurent polynomials, and corrects the ambient conjugate
variable to obtain a real derivation $\delta$ satisfying the necessary
inequalities.

\subsection{A cocycle with small coefficients}

The following lemma is useful in controlling sums occurring when applying
a cocycle to a long word in the free group.
\begin{lem}
\label{lem:sparse-clt}For each $m$, let $g_{m,1},\ldots,g_{m,N_{m}}$
belong to a discrete group $G_{m}$, with $N_{m}\to\infty$. Suppose
that for every fixed $p\geq1$, and all sufficiently large $m$, the
following hold:
\begin{enumerate}
\item If a word of length at most $2p$ in the letters $g_{m,j}^{\pm1}$
equals the identity, then the total degree of each letter must be
zero.
\item All but $O(N_{m}^{p-1})$ of the $p$-element subsets of letters are
free inside $G_{m}$.
\end{enumerate}
Let $\lambda_{m}$ be the left regular representation of $G_{m}$.
Then for any signs $\varepsilon_{m,j}\in\{-1,1\}$, the sums 
\[
Z_{m}=\frac{1}{\sqrt{N_{m}}}\sum_{j=1}^{N_{m}}\varepsilon_{m,j}\lambda_{m}(g_{m,j})
\]
converge in $*$-distribution to a standard circular element.
\end{lem}

\begin{proof}
Let
\[
H_{m}=\frac{1}{\sqrt{N_{m}}}\sum_{j=1}^{N_{m}}\varepsilon_{m,j}u_{j},
\]
where the $u_{j}$ are freely independent Haar unitaries.

Let $q(x,x^{*})=x^{\epsilon(1)}x^{\epsilon(2)}\cdots x^{\epsilon(r)}$
be a $*$-monomial of degree $r$, where $\epsilon:\{1,\dots,r\}\to\{\cdot,*\}$.
Then
\begin{align}
\tau(q(Z_{m},Z_{m}^{*})) & =N_{m}^{-\frac{r}{2}}\sum_{j_{1},\dots,j_{r}=1}^{N_{m}}\varepsilon_{m,j_{1}}\cdots\varepsilon_{m,j_{r}}\tau(\lambda_{m}(g_{m,j_{1}})^{\epsilon(1)}\cdots\lambda_{m}(g_{m,j_{r}})^{\epsilon(r)})\label{eq:expansionZm}\\
\tau(q(H_{m},H_{m}^{*})) & =N_{m}^{-\frac{r}{2}}\sum_{j_{1},\dots,j_{r}=1}^{N_{m}}\varepsilon_{m,j_{1}}\cdots\varepsilon_{m,j_{r}}\tau(u_{j_{1}}^{\epsilon(1)}\cdots u_{j_{r}}^{\epsilon(r)}).\label{eq:expansionHm}
\end{align}
Because of (1), for large enough $m$, the term $\tau(\lambda_{m}(g_{m,j_{1}})^{\epsilon(1)}\cdots\lambda_{m}(g_{m,j_{r}})^{\epsilon(r)})$
vanishes unless the total degree of each letter $g_{m,k}$ is zero.
The same condition is necessary for $\tau(u_{j_{1}}^{\epsilon(1)}\cdots u_{j_{r}}^{\epsilon(r)})$
to be nonzero, by freeness. Thus, for nonzero terms, it is possible
to pair up the integers $\{1,\dots,r\}$ in such a way that if $k$
is paired with $k'$ then $j_{k}=j_{k'}$ and $\epsilon(k)\neq\epsilon(k')$.
Thus it must be that $r$ is even (so $r=2p$ for some $p$); moreover,
there cannot be more than $p$ distinct letters. Thus we only need
to keep terms in the summations (\ref{eq:expansionZm}) and (\ref{eq:expansionHm})
in which the indices $j_{1},\dots,j_{2p}$ take on at most $p$ distinct
values.

There are at most $O(N_{m}^{p-1})$ terms in which the indices $j_{1},\dots,j_{2p}$
take on at most $p-1$ distinct values.

By (2), among terms in (\ref{eq:expansionZm}) in which the indices
$j_{1},\dots,j_{2p}$ take on exactly $p$ distinct values, all but
$O(N_{m}^{p-1})$ involve letters that freely generate a free subgroup,
and therefore equal the corresponding term in (\ref{eq:expansionHm}).

It follows that the difference between (\ref{eq:expansionZm}) and
(\ref{eq:expansionHm}) involves a summation over at most $O(N_{m}^{p-1})$
terms, each of which has absolute value at most one. Noting the prefactor
$N_{m}^{-r/2}=N_{m}^{-p}$, we see that for large enough $m$,
\[
\tau(q(Z_{m},Z_{m}^{*}))=\tau(q(H_{m},H_{m}^{*}))+O(N_{m}^{-1}).
\]

The usual free central limit theorem (see e.g. \cite[Lectures 11--12]{NicaSpeicher})
then implies that the limit law of $(H_{m},H_{m}^{*})$ (and thus
of $(Z_{m},Z_{m}^{*})$) is the circular law.
\end{proof}
\begin{lem}
\label{lem:cocycle}Let $G=\mathbb{F}_{2}=\langle a,b\rangle$, and
let $\lambda$ be its left regular representation on $\ell^{2}(G)$.
For every $\eta>0$, there are a finitely supported real vector $x\in\ell^{2}(G)$,
and $w\in G\setminus\{e\}$ such that the cocycle $c:G\to\ell^{2}(G)$
determined by 
\[
c(a)=(1-\lambda(a))x,\qquad c(b)=0
\]
satisfies 
\[
\|x\|_{2}<\eta,\qquad\|c(w)-\delta_{e}\|_{2}<\eta.
\]
\end{lem}

\begin{proof}
Let $x\in\ell^{2}(G)$ be an arbitrary real coefficient vector to
be chosen later, and let $c$ be the cocycle determined by $c(a)=(1-\lambda(a))x$
and $c(b)=0$.

For $m\geq1$, put 
\[
M=m+2,\qquad r_{j}=M^{j}-M^{j-1},\qquad w_{m}=\prod_{j=1}^{m}a^{r_{j}}b^{r_{j}}.
\]
Write 
\[
g_{j}=\prod_{i<j}a^{r_{i}}b^{r_{i}},\qquad h_{j}=g_{j}a^{r_{j}},\qquad R_{m}=\sum_{j=1}^{m}\bigl(\lambda(g_{j})-\lambda(h_{j})\bigr).
\]
Since $c$ is inner on $\langle a\rangle$, $c(a^{r_{j}})=(1-\lambda(a^{r_{j}}))x$.
Since $c(b)=0$, the cocycle identity implies that $c(w_{m})=R_{m}x$.

Let $\tau$ be the canonical trace on $L(G)$. We apply Lemma \ref{lem:sparse-clt}
to the $2m-1$ group elements 
\[
g_{2},\ldots,g_{m},\ h_{1},\ldots,h_{m}.
\]
For condition (1), note that the total degree of $a$ or $b$ in $g_{j}$
is $\sum_{i<j}r_{i}=M^{j-1}-1$. The total degree of $a$ in $h_{j}$
is $M^{j}-1$, and the total degree of $b$ in $h_{j}$ is $M^{j-1}-1$.
Thus, if a word in the letters $g_{2},\dots,g_{m},h_{1},\dots,h_{m}$
of length at most $2p$ reduces to the identity, and the total degree
of $g_{k}$ is $t_{k}$ and the total degree of $h_{k}$ is $s_{k}$
(with $t_{1}=0$), then
\[
\sum t_{k}(M^{k-1}-1)+\sum s_{k}(M^{k-1}-1)=0
\]
(because the total degree of $b$ must be zero) and
\[
\sum t_{k}(M^{k-1}-1)+\sum s_{k}(M^{k}-1)=0
\]
(because the total degree of $a$ must be zero). Subtracting these
we get
\[
(M-1)\sum s_{k}M^{k-1}=0.
\]
Since each $|s_{k}|\leq2p$, this forces $s_{k}=0$ once $M$ is large
enough. It follows that also $\sum t_{k}(M^{k-1}-1)=0$, which, because
$|t_{k}|\leq2p$, also forces $t_{k}=0$ for large enough $M$.

For condition (2), choose
\[
\begin{gathered}3\leq j_{1}<\cdots<j_{p}\leq m,\qquad j_{i+1}-j_{i}\geq2\quad(1\leq i<p),\\
v_{i}\in\{g_{j_{i}},h_{j_{i}}\}\qquad(1\leq i\leq p).
\end{gathered}
\]
Only $O_{p}(m^{p-1})$ selections fail these conditions. Write
\[
\begin{gathered}s_{2j-1}=a^{r_{j}},\qquad s_{2j}=b^{r_{j}},\qquad n_{0}=0,\\
v_{i}=s_{1}\cdots s_{n_{i}},\qquad n_{i}\in\{2j_{i}-2,2j_{i}-1\}.
\end{gathered}
\]
With $v_{0}=e$, put
\[
d_{i}=v_{i-1}^{-1}v_{i}=s_{n_{i-1}+1}\cdots s_{n_{i}},\qquad\ell(d_{i})=n_{i}-n_{i-1}\geq3,
\]
where $\ell$ counts the nonzero powers in the reduced alternating
expression in $a,b$. Since $r_{j}\ne r_{k}$ for $j\ne k$, complete
cancellation of adjacent powers in a product of the $d_{i}^{\pm1}$
can occur only between $d_{i}$ and $d_{i}^{-1}$. Thus, for $\varepsilon_{\nu}\in\{-1,1\}$
with
\[
(i_{\nu+1},\varepsilon_{\nu+1})\ne(i_{\nu},-\varepsilon_{\nu})\qquad(1\leq\nu<k),
\]
we have
\[
\ell\!\left(d_{i_{1}}^{\varepsilon_{1}}\cdots d_{i_{k}}^{\varepsilon_{k}}\right)\geq\sum_{\nu=1}^{k}\ell(d_{i_{\nu}})-(k-1)\geq2k+1>0.
\]
Hence $d_{i}$ are free inside $\mathbb{F}_{2}$, and so are the $v_{i}=d_{1}\cdots d_{i}$.

Consider now the signed sum
\[
\widetilde{R}_{m}=\sum_{j=2}^{m}\lambda(g_{j})-\sum_{j=1}^{m}\lambda(h_{j})=R_{m}-1.
\]
Lemma \ref{lem:sparse-clt} now shows that $(2m-1)^{-1/2}\widetilde{R}_{m}$
(and thus $(2m-1)^{-1/2}R_{m}$) converges in $*$-distribution to
the circular law. Consequently the spectral distributions $\nu_{m}$
of 
\[
A_{m}=\frac{R_{m}R_{m}^{*}}{2m}
\]
converge weakly to the free Poisson law
\[
d\nu(s)=\frac{1}{2\pi}\sqrt{\frac{4-s}{s}}\,\mathbf{1}_{(0,4)}(s)\,ds
\]
(see \cite[Example 15.15 and Proposition 12.13]{NicaSpeicher}). In
particular, $\nu(\{0\})=0$.

Now let 
\[
\widetilde{x}_{m}=R_{m}^{*}(1+R_{m}R_{m}^{*})^{-1}\delta_{e}.
\]
Since 
\[
R_{m}\widetilde{x}_{m}-\delta_{e}=-(1+R_{m}R_{m}^{*})^{-1}\delta_{e},
\]
the spectral theorem gives 
\[
\|\widetilde{x}_{m}\|_{2}^{2}+\|R_{m}\widetilde{x}_{m}-\delta_{e}\|_{2}^{2}=\tau\bigl((1+R_{m}R_{m}^{*})^{-1}\bigr).
\]
Thus, for every fixed $\varepsilon>0$,
\[
\begin{aligned}\tau\bigl((1+R_{m}R_{m}^{*})^{-1}\bigr) & =\int\frac{1}{1+2ms}\,d\nu_{m}(s)\\
 & \leq\nu_{m}([0,\varepsilon])+\frac{1}{1+2m\varepsilon}.
\end{aligned}
\]

Since $\nu(\{0\})=0$, we may choose $\varepsilon>0$ such that $\nu([0,\varepsilon])<\eta^{2}/8$.
Weak convergence gives
\[
\limsup_{m\to\infty}\nu_{m}([0,\varepsilon])\leq\nu([0,\varepsilon])<\frac{\eta^{2}}{8}.
\]
We can therefore choose $m$ sufficiently large that both $\nu_{m}([0,\varepsilon])<\eta^{2}/8$
and $(1+2m\varepsilon)^{-1}<\eta^{2}/8$. The preceding identity and
bound then give
\[
\|\widetilde{x}_{m}\|_{2}^{2}+\|R_{m}\widetilde{x}_{m}-\delta_{e}\|_{2}^{2}<\frac{\eta^{2}}{4},
\]
so each of the two norms is less than $\eta/2$.

The vector $\widetilde{x}_{m}$ is real. Approximate $\widetilde{x}_{m}$
by a finitely supported real vector $x$ with 
\[
\|x-\widetilde{x}_{m}\|_{2}<\frac{\eta}{2(1+\|R_{m}\|)}.
\]
Then $\|x\|_{2}<\eta$ and $\|R_{m}x-\delta_{e}\|_{2}<\eta$. Since
$w_{m}$ is a nonempty positive word and $c(w_{m})=R_{m}x$, this
proves the lemma.
\end{proof}

\subsection{Norm estimates for the derivation}

We construct a scalar polynomial model, lift it to $B_{0}$, and correct
its kernel to retain the same conjugate variable exactly. The resulting
estimates give the derivation required by Proposition \ref{prop:absorption}.

We first establish two properties of $\|\cdot\|_{\mathrm{rc}}$. Let
$(u_{h})_{h\in\mathbb{F}_{2}}$ be the group unitaries, fix $g\in\mathbb{F}_{2}\setminus\{e\}$,
and let $f\in\ell^{2}(\mathbb{F}_{2})$ be finitely supported. Then
\[
\|\Delta(f)\#(u_{g}^{r}\otimes(u_{g}^{s})^{\mathrm{op}})\|_{\mathrm{rc}}=2\|f\|_{2}\qquad(r,s\in\mathbb{Z}).
\]
This is immediate from orthonormality of $(u_{h}u_{g}^{r})_{h}$ and
$(u_{g}^{s}u_{h}^{*})_{h}$ and the definition of $\Vert\cdot\Vert_{rc}$.
In particular, if $K(u_{g})=\sum_{r,s}\kappa_{rs}u_{g}^{r}\otimes(u_{g}^{s})^{\mathrm{op}}$,
then
\[
\|\Delta(f)\#K(u_{g})\|_{\mathrm{rc}}\leq2\left(\sum_{r,s}|\kappa_{rs}|\right)\|f\|_{2}.
\]
Moreover, for $H\in C(\mathbb{T}^{2})$ and any Haar unitary $V$
(identified with the complex coordinate function on $\mathbb{T}$)
\[
\begin{aligned}\|H(V\otimes1,1\otimes V^{\mathrm{op}})\|_{\mathrm{rc}}={} & \sup_{\zeta_{2}}\left(\int|H(\zeta_{1},\zeta_{2})|^{2}\,dm(\zeta_{1})\right)^{1/2}\\
 & +\sup_{\zeta_{1}}\left(\int|H(\zeta_{1},\zeta_{2})|^{2}\,dm(\zeta_{2})\right)^{1/2}\leq2\Vert H\Vert_{\infty}
\end{aligned}
\]
where $m$ is normalized Haar measure. This holds in particular for
$V=u_{g}$ whenever $g\in\mathbb{F}_{2}\setminus\{e\}$.
\begin{lem}
\label{lem:scalar-conjugate-variable}Let $V$ be a Haar unitary.
Identify $C^{*}(V)$ with $C(\mathbb{T})$ by identifying $V$ with
the complex coordinate, and regard tensors in $C^{*}(V)\otimes C^{*}(V)^{\mathrm{op}}$
as functions on $\mathbb{T}^{2}$. Denote by $m$ the Haar measure.

For every $\epsilon_{0}>0$, there is a real-valued antisymmetric
Laurent polynomial $q$ such that the inner derivation $d$ on $\mathop{{*}\text{-}\mathrm{Alg}}\nolimits(V)$
given by $d(g)=gq-qg$ and the selfadjoint element $p\in\mathop{{*}\text{-}\mathrm{Alg}}\nolimits(V)$
given by $p(\zeta)=2\int q(\zeta,\omega)\,dm(\omega)$ satisfy 
\[
p=d^{*}(1\otimes1),\qquad\|d(p)-1\otimes1\|_{\mathrm{rc}}<\epsilon_{0}.
\]
\end{lem}

\begin{proof}
Given a continuous real-valued antisymmetric $q\in C(\mathbb{T}^{2})$,
let $T_{q}$ be its integral operator on $L^{2}(\mathbb{T},m)$. Under
the standard identification with Hilbert-{}-Schmidt operators, the
inner derivation $dg=gq-qg$ becomes $[M_{g},T_{q}]$, where $M_{g}$
is multiplication by $g$. Since $T_{q}$ is real and $T_{q}^{*}=-T_{q}$,
the commutator formula of \cite[Proposition 2.6]{ShlyakhtenkoPartialConjugates}
gives 
\[
d^{*}(1\otimes1)=2T_{q}1=2\int q(\,\cdot\,,\omega)\,dm(\omega).
\]

Let $f:\mathbb{T}\to[-2,2]$ be continuous and such that $f_{*}m=\mu$
is the semicircular measure. Write $x=f(\zeta)$, $y=f(\omega)$.
For $t>0$ let 
\[
\begin{aligned}q_{t}(x,y) & =\frac{x-y}{(x-y)^{2}+t^{2}},\\
h_{t}(x) & =x-2\int q_{t}(x,y)\,d\mu(y).
\end{aligned}
\]
The semicircle Cauchy transform \cite[Lemma 2.21]{NicaSpeicher} gives
\[
h_{t}(x)=\operatorname{Re}\sqrt{(x+it)^{2}-4},\qquad\|h_{t}\|_{\infty}\leq\sqrt{2t},
\]
using the branch $\sqrt{z^{2}-4}\sim z$ at infinity. Also, $\int h_{t}\,d\mu=0$
by centering and antisymmetry.

For a kernel $q$ and a real function $r$, write
\[
q[r](x,y)=q(x,y)+\frac{r(x)-r(y)}{2}.
\]
If $\int r\,d\mu=0$, this correction changes $2\int q(x,y)\,d\mu(y)$
by $r(x)$. Take $\widetilde{q}_{t}=q_{t}[h_{t}]$, and let $\widetilde{d}_{t}$
be the derivation associated with the kernel $\widetilde{q}_{t}(f(\zeta),f(\omega))$.
Then $\widetilde{d}_{t}^{*}(1\otimes1)=f$. Moreover, 
\[
\widetilde{d}_{t}f-1\otimes1=-\frac{t^{2}}{(x-y)^{2}+t^{2}}+\frac{(x-y)(h_{t}(x)-h_{t}(y))}{2}.
\]
Since the semicircle density is at most $1/\pi$,
\[
\left\Vert \frac{t^{2}}{(x-y)^{2}+t^{2}}\right\Vert _{rc}^{2}\leq\frac{4}{\pi}\int\frac{t^{4}}{(u^{2}+t^{2})^{2}}du=2t
\]
and
\[
\left\Vert \frac{(x-y)(h_{t}(x)-h_{t}(y))}{2}\right\Vert _{rc}\leq4\Vert x\Vert_{\infty}\left\Vert h_{t}\right\Vert _{\infty}\leq8\sqrt{2t}
\]
so that
\[
\|\widetilde{d}_{t}f-1\otimes1\|_{\mathrm{rc}}\leq9\sqrt{2t}.
\]

Fix $t$ with $9\sqrt{2t}<\epsilon_{0}/2$. Approximate the continuous
kernel $\widetilde{q}_{t}(f(\zeta),f(\omega))$ uniformly by real-valued
antisymmetric Laurent polynomials $q_{n}$. Put 
\[
p_{n}(\zeta)=2\int q_{n}(\zeta,\omega)\,dm(\omega),\qquad d_{n}(g)(\zeta,\omega)=(g(\zeta)-g(\omega))q_{n}(\zeta,\omega).
\]
Then $p_{n}=d_{n}^{*}(1\otimes1)$ and, as $n\to\infty$, 
\[
p_{n}\longrightarrow f,\qquad d_{n}(p_{n})\longrightarrow\widetilde{d}_{t}f
\]
uniformly. Since $\|H\|_{\mathrm{rc}}\leq2\|H\|_{\infty}$, the desired
estimate holds for $p=p_{n}$, $q=q_{n}$ and $d=d_{n}$ when $n$
is sufficiently large.
\end{proof}
For Laurent polynomials in one or two variables, define
\[
\left\Vert \sum_{j}f_{j}\zeta^{j}\right\Vert _{\mathrm{W}}=\sum_{j}|f_{j}|,\qquad\left\Vert \sum_{j,k}f_{j,k}\zeta^{j}\omega^{k}\right\Vert _{\mathrm{W}}=\sum_{j,k}|f_{j,k}|.
\]

We next use this computation to produce elements in $B_{0}=\mathbb{C}G$.
Fix $\epsilon_{0}>0$ and choose $q,p$ as in Lemma \ref{lem:scalar-conjugate-variable}.
Write $p(\zeta)=\sum_{j}p_{j}\zeta^{j},$ $q(\zeta,\omega)=\sum q_{j,k}\zeta^{j}\omega^{k}$,
and put
\[
A_{q}=\|q\|_{\mathrm{W}}=\sum_{j,k}|q_{j,k}|,\qquad L_{p}=\sum_{j}|j|\,|p_{j}|,\qquad C_{q}=8A_{q}(L_{p}+2\|p\|_{\infty}),
\]

\begin{prop}
\label{prop:lifted-conjugate-variable}Let $x\in B_{0}$, $w\in G\setminus\{e\}$
and $c:G\to\ell^{2}(G)$ be as in Lemma \ref{lem:cocycle}. Put
\[
\rho=\|x\|_{2},\qquad\theta=\|c(w)-\delta_{e}\|_{2},
\]
and assume that $\theta<1/4$.

For a real-valued antisymmetric Laurent polynomial $\hat{q}(\zeta,\omega)=\sum\hat{q}_{j,k}\zeta^{j}\omega^{k}$,
let $\hat{K}=(\zeta-\omega)\hat{q}(\zeta,\omega)\omega^{-1}$, $Q=\Delta(x)\#\hat{K}(u_{w})$,
and consider the derivation $\delta:B_{0}\to B_{0}\otimes B_{0}^{\mathrm{op}}$
defined by
\[
\delta(u_{a})=Qu_{a}-u_{a}Q,\qquad\delta(u_{b})=0.
\]
Then we can choose $\hat{q}$ so that $\Vert\hat{q}-q\Vert_{W}\leq4A_{q}\theta$,
$Q^{\dagger}=-Q$, $\delta$ is real, $\xi=\delta^{*}(1\otimes1)=p(u_{w})$,
and the estimates
\[
\|Q\|_{\mathrm{rc}}\leq8A_{q}\rho,\qquad\|\delta(u_{a})\|_{\mathrm{rc}}\leq16A_{q}\rho,\qquad\|\delta(\xi)-1\otimes1\|_{\mathrm{rc}}\leq\epsilon_{0}+C_{q}\theta
\]
hold.
\end{prop}

\begin{proof}
For any real antisymmetric polynomial $q'$, put $K_{q'}(\zeta,\omega)=(\zeta-\omega)q'(\zeta,\omega)\omega^{-1}$.
Since $K_{q'}^{\dagger}=-K_{q'}$ and $x$ has real coefficients,
$Q_{q'}=\Delta(x)\#K_{q'}(u_{w})$ satisfies $Q_{q'}^{\dagger}=-Q_{q'}$.
Define $\delta_{q'}$ on $B_{0}$ by $\delta_{q'}(u_{a})=[Q_{q'},u_{a}]$
and $\delta_{q'}(u_{b})=0$. The cocycle identity then gives
\[
\delta_{q'}(u_{g})=(\Delta(c(g))\#K_{q'}(u_{w}))u_{g}\qquad(g\in G).
\]

Set
\[
\gamma_{n}=c(w)(w^{n}),\qquad\sigma_{n}=\frac{\gamma_{n}+\gamma_{-n}}{2},\qquad F(\zeta,\omega)=\sum_{n}\gamma_{n}\zeta^{n}\omega^{-n}.
\]
These coefficients are real and finitely supported. For a Laurent
polynomial $g$ the chain rule and definition of $\delta_{q'}$ give:
\[
(\tau\otimes\tau)(\delta_{q'}(g(u_{w})))=\iint(g(\zeta)-g(\omega))F(\zeta,\omega)q'(\zeta,\omega)\,dm(\zeta)\,dm(\omega).
\]
Since $F(\omega,\zeta)=\overline{F(\zeta,\omega)}$ and $q'$ is antisymmetric,
the real centered Laurent polynomial
\[
\xi_{q'}(\zeta)=2\int q'(\zeta,\omega)\operatorname{Re}F(\zeta,\omega)\,dm(\omega),\qquad\operatorname{Re}F=\sum_{n}\sigma_{n}\zeta^{n}\omega^{-n}
\]
satisfies
\[
(\tau\otimes\tau)(\delta_{q'}z)=\tau(\xi_{q'}(u_{w})z)\qquad\bigl(\forall z\in\mathop{{*}\text{-}\mathrm{Alg}}\nolimits(u_{w})\bigr).
\]

This identity holds on all of $B_{0}$. Indeed, write $K_{q'}(\zeta,\omega)=\sum_{r,s}\kappa_{rs}\zeta^{r}\omega^{s}$.
Each tensor in the expansion of $\delta_{q'}(u_{g})$ has the form
\[
u_{h}u_{w}^{r}\otimes(u_{w}^{s}u_{h}^{*}u_{g})^{\mathrm{op}}.
\]
Both factors have nonzero trace only if $h=w^{-r}$ and $g=w^{-(r+s)}$.
Thus, for $g\notin\langle w\rangle$, both sides of the preceding
identity vanish. By linearity the identity holds for every $z\in\mathbb{C}G$,
proving that $1\otimes1\in\operatorname{Dom}(\delta_{q'}^{*})$ and
\[
\delta_{q'}^{*}(1\otimes1)=\xi_{q'}(u_{w})
\]
where the adjoint of $\delta_{q'}$ is computed by regarding $\delta_{q'}$
as an unbounded operator $\delta_{q'}:B_{0}\to L^{2}(L(G)\otimes L(G)^{op})$.

Apply the correction $q[r]$ from Lemma \ref{lem:scalar-conjugate-variable},
now with the weight $\operatorname{Re}F$. For a centered real polynomial
$r(\zeta)=\sum_{n\ne0}r_{n}\zeta^{n}$, the preceding formula gives
\[
\xi_{q[r]}(\zeta)-\xi_{q}(\zeta)=\sum_{n\ne0}(\gamma_{0}-\sigma_{n})r_{n}\zeta^{n}.
\]
The bounds $|\gamma_{0}-1|\leq\theta$ and $|\sigma_{n}|\leq\theta$
for $n\ne0$ imply
\[
\|p-\xi_{q}\|_{\mathrm{W}}\leq2A_{q}\theta,\qquad\gamma_{0}-\sigma_{n}\geq1-2\theta>1/2\quad(n\ne0).
\]
Writing $p-\xi_{q}=\sum_{n\ne0}\beta_{n}\zeta^{n}$, set
\[
r(\zeta)=\sum_{n\ne0}\frac{\beta_{n}}{\gamma_{0}-\sigma_{n}}\zeta^{n},\qquad\hat{q}=q[r].
\]
Then $r$ is real and centered, $\|r\|_{\mathrm{W}}\leq4A_{q}\theta$,
and $\xi_{\hat{q}}=p$. Take $\delta=\delta_{\hat{q}}$.

Since $\theta<1/4$, we have $\|\hat{q}\|_{\mathrm{W}}\leq2A_{q}$
and $\|\hat{K}\|_{\mathrm{W}}\leq4A_{q}$. The coefficient estimate
above and the bimodule inequality therefore give
\[
\|Q\|_{\mathrm{rc}}\leq2\rho\|\hat{K}\|_{\mathrm{W}}\leq8A_{q}\rho,\qquad\|\delta(u_{a})\|_{\mathrm{rc}}\leq2\|Q\|_{\mathrm{rc}}.
\]
Let $\hat{d}(g)=(g(\zeta)-g(\omega))\hat{q}(\zeta,\omega)$ be the
derivation defined on $\mathop{{*}\text{-}\mathrm{Alg}}\nolimits(u_{w})$.
The identity
\[
\delta(u_{w})-\hat{d}(u_{w})=(\Delta(c(w)-\delta_{e})\#\hat{K}(u_{w}))u_{w}
\]
and the Leibniz rule give
\[
\|\delta(p(u_{w}))-\hat{d}(p(u_{w}))\|_{\mathrm{rc}}\leq2\theta\|\hat{K}\|_{\mathrm{W}}L_{p}\leq8A_{q}\theta L_{p}.
\]

Finally,
\[
((\hat{d}-d)p)(\zeta,\omega)=\frac{(p(\zeta)-p(\omega))(r(\zeta)-r(\omega))}{2}.
\]
Using $\|H\|_{\mathrm{rc}}\leq2\|H\|_{\infty}$ and Lemma \ref{lem:scalar-conjugate-variable},
we obtain
\[
\begin{aligned}\|\delta(\xi)-1\otimes1\|_{\mathrm{rc}} & \leq\epsilon_{0}+8A_{q}\theta L_{p}+4\|p\|_{\infty}\|r\|_{\mathrm{W}}\\
 & \leq\epsilon_{0}+C_{q}\theta.
\end{aligned}
\]
\end{proof}
\begin{cor}
\label{cor:conjugate-variable-data}For every $\eta,\epsilon,\kappa>0$,
there is a real derivation $\delta:\mathbb{C}G\to\mathbb{C}G\otimes(\mathbb{C}G)^{\mathrm{op}}$
such that
\[
\delta(u_{b})=0,\qquad\xi=\delta^{*}(1\otimes1)\in\mathbb{C}G,
\]
\[
\|\delta(u_{a})\|_{\mathrm{rc}}<\eta,\qquad\|\delta(\xi)-1\otimes1\|_{\mathrm{rc}}<\epsilon.
\]
Moreover, $\delta$ is inner on $\mathop{{*}\text{-}\mathrm{Alg}}\nolimits(u_{a})$,
given by $\delta(y)=[Q,y]$, $y\in\mathop{{*}\text{-}\mathrm{Alg}}\nolimits(u_{a})$,
with $Q$ satisfying $Q^{\dagger}=-Q$ and $\|Q\|_{\mathrm{rc}}<\kappa$.
\end{cor}

\begin{proof}
First apply Lemma \ref{lem:scalar-conjugate-variable} with $\epsilon_{0}<\epsilon/2$,
fixing $q,p$ and hence $A_{q},C_{q}$. Choose $0<\eta_{0}<1/4$ so
that
\[
16A_{q}\eta_{0}<\eta,\qquad8A_{q}\eta_{0}<\kappa,\qquad C_{q}\eta_{0}<\epsilon/2.
\]
Apply Lemma \ref{lem:cocycle} with tolerance $\eta_{0}$. Then $\rho,\theta<\eta_{0}$,
and Proposition \ref{prop:lifted-conjugate-variable} gives the required
derivation, with $\xi=p(u_{w})$ exactly and with the asserted small
inner implementation. All coefficients and the conjugate variable
lie in the group algebra.
\end{proof}

\section{Semicircular absorption}\label{sec:Semicircular-absorption}
\begin{thm}
\label{thm:free-absorption}Let $2\leq n<\infty$, let $u_{1},\ldots,u_{n}$
be the canonical free Haar generators of $B_{n}=C_{\mathrm{red}}^{*}(\mathbb{F}_{n})$,
and put $A_{n}=C^{*}(u_{2},\ldots,u_{n})$. If $S$ is standard semicircular
and free from $B_{n}$, then
\[
B_{n}\subset C_{n}=B_{n}*_{\mathrm{red}}C^{*}(S)
\]
has CAAP relative to $A_{n}$, with small unitary implementation on
$C^{*}(u_{1})$.

Equivalently, there are surjective trace-preserving isomorphisms $\Phi_{j}:B_{n}\to C_{n}$
fixing $A_{n}$ pointwise and converging pointwise in norm to inclusion,
with
\[
\Phi_{j}(u_{1})=U_{j}u_{1}U_{j}^{*},\qquad U_{j}\in\mathcal{U}(C_{n}),\qquad\|U_{j}-1\|\longrightarrow0.
\]

In particular,
\[
B_{n}\cong B_{n}*_{\mathrm{red}}C^{*}(S).
\]

Passing to the GNS representation, the inclusion
\[
L(\mathbb{F}_{n})\subset L(\mathbb{F}_{n})*L(\mathbb{Z})\cong L(\mathbb{F}_{n+1})
\]
has AAP relative to $W^{*}(u_{2},\ldots,u_{n})$, with small unitary
implementation on $W^{*}(u_{1})$. In particular, $L(\mathbb{F}_{n})\cong L(\mathbb{F}_{n+1})$.
\end{thm}

\begin{proof}
For $n=2$, apply Corollary \ref{cor:conjugate-variable-data} with
$u_{a}=u_{1}$ and $u_{b}=u_{2}$. Its derivations satisfy the hypotheses
of Proposition \ref{prop:absorption} on $B_{0}=\mathbb{C}\mathbb{F}_{2}$,
with $A=\mathop{{*}\text{-}\mathrm{Alg}}\nolimits(u_{2})$ and $X=u_{1}$,
and the bound for the inner restriction to $\mathop{{*}\text{-}\mathrm{Alg}}\nolimits(u_{1})$.
The proposition gives CAAP relative to $C^{*}(u_{2})$, with small
unitary implementation on $C^{*}(u_{1})$. By Proposition \ref{prop:caap},
we may choose surjective trace-preserving isomorphisms
\[
\Psi_{j}:B_{2}\longrightarrow B_{2}*_{\mathrm{red}}C^{*}(S)
\]
fixing $u_{2}$ and converging pointwise in norm to inclusion, and
such that $\Psi_{j}(u_{1})=U_{j}u_{1}U_{j}^{*}$ and $\|U_{j}-1\|\to0$.

For $n>2$, put $D=C^{*}(u_{3},\ldots,u_{n})$, so $B_{n}=B_{2}*_{\mathrm{red}}D$.
We therefore obtain surjective trace-preserving isomorphisms
\[
\begin{aligned}\Phi_{j}=\Psi_{j}*\mathrm{id}_{D}:\ B_{n} & \longrightarrow(B_{2}*_{\mathrm{red}}C^{*}(S))*_{\mathrm{red}}D\\
 & \cong(B_{2}*_{\mathrm{red}}D)*_{\mathrm{red}}C^{*}(S)=C_{n}.
\end{aligned}
\]
The maps $\Phi_{j}$ fix $u_{2}$ and $D$, hence fix $A_{n}$, and
retain $\Phi_{j}(u_{1})=U_{j}u_{1}U_{j}^{*}$ with the same unitaries,
and converge pointwise to the inclusion $B_{n}\subset C_{n}$. The
rest of the statements are straightforward. For the von Neumann algebra
application, it only remains to pass to the von Neumann algebras,
noting that $W^{*}(S)\cong L(\mathbb{Z})$ since the semicircular
measure is diffuse.
\end{proof}
\begin{cor}
\label{cor:countable-absorption}Let $B_{n},A_{n}$ be as in Theorem
\ref{thm:free-absorption}, and let $(S_{j})_{j\in J}$ be a nonempty
finite or countable free standard semicircular family, free from $B_{n}$.
The inclusion
\[
B_{n}\subset C^{*}(B_{n},S_{j}:j\in J)
\]
has CAAP relative to $A_{n}$, with small unitary implementation on
$C^{*}(u_{1})$. Consequently,
\[
B_{n}\cong B_{n}*_{\mathrm{red}}C^{*}(S_{j}:j\in J).
\]
The generated von Neumann inclusion has AAP relative to $W^{*}(A_{n})$,
with small unitary implementation on $W^{*}(u_{1})$. If $J$ is infinite,
this gives $L(\mathbb{F}_{n})\cong L(\mathbb{F}_{\infty})$.
\end{cor}

\begin{proof}
Enumerate the family and put $D_{0}=B_{n}$, $D_{k}=C^{*}(B_{n},S_{1},\ldots,S_{k})$,
and $E=C^{*}(u_{1})$. Theorem \ref{thm:free-absorption} gives relative
CAAP for $D_{0}\subset D_{1}$ with small unitary implementation on
$E$. Inductively, suppose that $B_{n}\subset D_{k}$ has this strengthened
property. Choose a surjective trace-preserving isomorphism $\Theta:B_{n}\to D_{k}$
fixing $A_{n}$ and satisfying $\Theta|_{E}=\operatorname{Ad}T|_{E}$
for a unitary $T\in D_{k}$. The isomorphism $\Theta*\mathrm{id}_{C^{*}(S_{k+1})}$
carries the inclusion $B_{n}\subset B_{n}*_{\mathrm{red}}C^{*}(S_{k+1})$
to $D_{k}\subset D_{k+1}$. Thus the latter has relative CAAP with
small unitary implementation on $\Theta(E)=TET^{*}$. Remark \ref{rem:implementation}
gives the same property on $E$, and Proposition \ref{prop:permanence}(1)
gives it for $B_{n}\subset D_{k+1}$. This proves the assertion for
finite families.

For an infinite family, put $D_{\infty}=\overline{\bigcup_{k}D_{k}}^{\|\cdot\|}$.
Proposition \ref{prop:permanence}(2) gives CAAP for $B_{n}\subset D_{\infty}$
relative to $A_{n}$, with small unitary implementation on $E$. Proposition
\ref{prop:caap} gives the stated $C^{*}$-isomorphism.
\end{proof}
Theorem \ref{thm:free-absorption} proves Examples \ref{exa:MainExamples}
(a) and (b), with the stated relative approximation properties. Example
(c) follows directly from \cite{ShlyakhtenkoFuchsian}.

\section{Further free-product consequences}\label{sec:free-product-consequences}

The small unitary implementation obtained in Corollary \ref{cor:countable-absorption}
is used below.
\begin{lem}
\label{lem:small-conjugation}Let $u,v$ be free Haar unitaries, and
let $(S_{j})_{j\in J}$ be a nonempty finite or countable free standard
semicircular family, free from them. Put $R=C^{*}(u,v,S_{j}:j\in J)$.
For every $\varepsilon>0$ there is a unitary $U\in R$ such that
\[
\|U-1\|<\varepsilon,\qquad R=C^{*}(UuU^{*},v),
\]
and $UuU^{*},v$ are free Haar unitaries.
\end{lem}

\begin{proof}
Apply Corollary \ref{cor:countable-absorption} with $n=2$, $u_{1}=u$
and $u_{2}=v$. By the strengthened form of Proposition \ref{prop:caap},
there are a surjective trace-preserving isomorphism $\Phi:C^{*}(u,v)\to R$
fixing $v$ and a unitary $U\in R$ such that
\[
\Phi(u)=UuU^{*},\qquad\|U-1\|<\varepsilon.
\]
Surjectivity gives $R=C^{*}(UuU^{*},v)$, and trace preservation gives
the free Haar law of $UuU^{*},v$.
\end{proof}
\begin{thm}
\label{thm:freeProductsAbsorb}Let $(A_{i},\tau_{i})$, $i=1,2$ be
separable $C^{*}$-algebras with faithful tracial states, and suppose
that there exist Haar unitaries $u\in A_{1}$ and $v\in A_{2}$. Let
$J$ be a finite or countable set, and let $S_{j}:j\in J$ be a free
semicircular family free from $A_{1}*_{\textrm{red}}A_{2}$.

Then for every $\epsilon>0$ there exists a unitary $U\in C^{*}(u,v,S_{j}:j\in J)$,
$\Vert U-1\Vert<\epsilon$, and so that $C^{*}(A_{1},A_{2},S_{j}:j\in J)=C^{*}(UA_{1}U^{*},A_{2})$
and moreover $UA_{1}U^{*}$ and $A_{2}$ are freely independent. Thus
$A_{1}*_{\textrm{red}}A_{2}\subset A_{1}*_{\textrm{red}}A_{2}*_{\textrm{red}}C^{*}(S_{j}:j\in J)$
has CAAP.

In particular, $A_{1}*_{\textrm{red}}A_{2}\cong A_{1}*_{\textrm{red}}A_{2}*_{\textrm{red}}C^{*}(S_{j}:j\in J)$.
\end{thm}

\begin{proof}
Clearly, we can assume that $J$ is nonempty. Let $Q=C^{*}(S_{j}:j\in J)$.
For every $\epsilon>0$, Lemma \ref{lem:small-conjugation} yields
a unitary $U\in C^{*}(u,v,S_{j}:j\in J)$ so that $\Vert U-1\Vert<\epsilon$,
and so that the map $\Phi:u\mapsto UuU^{*},\ v\mapsto v$ extends
to a $*$-isomorphism of $B=C^{*}(u,v)$ with $C^{*}(u,v,S_{j}:j\in J)=B*_{\textrm{red}}Q$,
and so that $UuU^{*}$ and $v$ are freely independent.

Let $\mathcal{A}_{j}=W^{*}(A_{j})$, $\mathcal{Q}=W^{*}(S_{j}:j\in J)$
in the GNS representation corresponding to tracial states. We repeatedly
use the identity $M*_{B}(B*N)=M*N$ valid for amalgamated free products
of von Neumann algebras taken with respect to trace-preserving conditional
expectations.

Let now $\mathcal{A}=\mathcal{A}_{1}*\mathcal{A}_{2}$. Then
\begin{align*}
\mathcal{A}*\mathcal{Q} & =\mathcal{A}_{1}*\mathcal{Q}*\mathcal{A}_{2}=W^{*}(\mathcal{A}_{1},u,v,Q)*_{W^{*}(v)}\mathcal{A}_{2}\\
 & =(\mathcal{A}_{1}*_{W^{*}(u)}W^{*}(u,v,Q))*_{W^{*}(v)}\mathcal{A}_{2}\\
 & =(U\mathcal{A}_{1}U^{*}*_{W^{*}(UuU^{*})}W^{*}(u,v,Q))*_{W^{*}(v)}\mathcal{A}_{2}
\end{align*}
because $U\in W^{*}(u,v,Q)$. Since $C^{*}(UuU^{*},v)=C^{*}(u,v,Q)$
also $W^{*}(UuU^{*},v)=W^{*}(u,v,Q)$; moreover $UuU^{*}$ and $v$
are free, so $W^{*}(u,v,Q)=W^{*}(UuU^{*})*W^{*}(v)$. Thus
\begin{align*}
\mathcal{A}*\mathcal{Q} & =(U\mathcal{A}_{1}U^{*}*_{W^{*}(UuU^{*})}W^{*}(u,v,Q))*_{W^{*}(v)}\mathcal{A}_{2}\\
 & =(U\mathcal{A}_{1}U^{*}*_{W^{*}(UuU^{*})}[W^{*}(UuU^{*})*W^{*}(v)])*_{W^{*}(v)}\mathcal{A}_{2}\\
 & =(U\mathcal{A}_{1}U^{*}*W^{*}(v))*_{W^{*}(v)}\mathcal{A}_{2}=W^{*}(U\mathcal{A}_{1}U^{*})*\mathcal{A}_{2}.
\end{align*}

It follows that $U\mathcal{A}_{1}U^{*}$ and $\mathcal{A}_{2}$ are
free in $\mathcal{A}*\mathcal{Q}$. But then $UA_{1}U^{*}$ and $A_{2}$
are free in $A_{1}*_{\textrm{red}}A_{2}*_{\textrm{red}}Q$. It follows
that the maps
\begin{align*}
x & \mapsto UxU^{*},\qquad x\in A_{1},\\
y & \mapsto y,\qquad y\in A_{2}
\end{align*}
extend to a trace-preserving $C^{*}$-monomorphism $\Psi_{\epsilon}$
of $A_{1}*_{\textrm{red}}A_{2}$ into $A_{1}*_{\textrm{red}}A_{2}*_{\textrm{red}}Q$.
Since $\{u,v\}\cup Q\subset C^{*}(UuU^{*},v)$, it follows that the
image of this monomorphism contains both $Q$ and $U$, and thus $A_{1}$
and $A_{2}$, and so is surjective; thus it gives an isomorphism of
$A_{1}*_{\textrm{red}}A_{2}*_{\textrm{red}}Q$ with $A_{1}*_{\textrm{red}}A_{2}$.

Taking $\epsilon\to0$, we see that $\Psi_{\epsilon}$ converges to
the inclusion $A_{1}*_{\textrm{red}}A_{2}\subset A_{1}*_{\textrm{red}}A_{2}*_{\textrm{red}}Q$. 
\end{proof}
\begin{cor}
Let $J\subset\mathbb{N}$ be a finite or countable set with $|J|\geq2$,
and let $S_{j}:j\in J$ be a free semicircular system. Then 
\[
C^{*}(S_{j}:j\in J)\cong C^{*}(S_{1},S_{2}).
\]
In particular, the isomorphism class of $Q=C^{*}(S_{j}:j\in J)$ does
not depend on the cardinality of $J$.
\end{cor}

\begin{proof}
Write, after reindexing, $J=J'\cup\{1,2\}$. This is Theorem \ref{thm:freeProductsAbsorb}
applied to $A_{1}=C^{*}(S_{1})$, $A_{2}=C^{*}(S_{2})$ and the semicircular
family $\{S_{j}:j\in J'\}$.
\end{proof}
\begin{cor}
Let $G$ and $H$ be two infinite discrete groups, $J$ be a finite
or countable set, and $S_{j}:j\in J$ a free semicircular family.
Let $Q=C^{*}(S_{j}:j\in J)$. Then 
\[
C_{\mathrm{red}}^{*}(G*H)*_{\mathrm{red}}Q\cong C_{\mathrm{red}}^{*}(G*H).
\]
\end{cor}

\begin{proof}
If $G$ is an infinite discrete group, then $C_{\textrm{red}}^{*}(G)$
contains a Haar unitary \cite[Proposition 8.1]{ThielHaarUnitaries}. 
\end{proof}
\begin{rem}
More general free products can sometimes be regrouped into free products
that contain Haar unitaries. For example, if four algebras $A_{1},\dots,A_{4}$
contain centered unitaries $z_{i}\in A_{i}$, then $z_{1}z_{2}$ and
$z_{3}z_{4}$ are free Haar unitaries. Thus in particular 
\[
\left(*_{\mathrm{red},\,i=1}^{4}(M_{n_{i}}(\mathbb{C}),\operatorname{tr}_{n_{i}})\right)*_{\mathrm{red}}Q\cong*_{\mathrm{red},\,i=1}^{4}(M_{n_{i}}(\mathbb{C}),\operatorname{tr}_{n_{i}}),\qquad n_{i}\geq2,
\]
where $Q$ as before is generated by a finite or countable free semicircular
family. 
\end{rem}

\begin{rem}
\label{rem:hirshberg-phillips-1}Let $I=(C([0,1]),\lambda)$, with
Lebesgue measure, and let $D$ be the countably infinite reduced free
product of $I$ with itself. Then $D\cong Q=C^{*}(S_{1},S_{2},\dots)$,
where $S_{1},S_{2},\dots$ are free semicircular variables.

Let $\mathcal{Z}$ be the Jiang--Su algebra. By a result of Hirshberg--Phillips
\cite[Theorem 3.1]{HirshbergPhillips},
\[
D*_{\textrm{red}}\mathcal{Z}\cong D
\]
(trace-preserving isomorphism). Since $\mathcal{Z}$ contains a Haar
unitary, we can apply Theorem \ref{thm:freeProductsAbsorb} to conclude
that
\begin{align*}
\mathcal{Z}*_{\textrm{red}}\mathcal{Z} & \cong\mathcal{Z}*_{\textrm{red}}\mathcal{Z}*_{\textrm{red}}D\\
 & \cong\mathcal{Z}*_{\textrm{red}}D*_{\textrm{red}}\mathcal{Z}*_{\textrm{red}}D\\
 & \cong D*_{\textrm{red}}D\cong D.
\end{align*}
This answers \cite[Question 4.5]{HirshbergPhillips}.

In particular, reduced free products of $n\geq2$ (or countably infinitely
many) copies of $\mathcal{Z}$ are isomorphic to $\mathcal{Z}*_{\textrm{red}}\mathcal{Z}\cong I*_{\textrm{red}}I$.
\end{rem}

\begin{rem}
\label{rem:nccw-absorption}More generally, let $E$ be a unital countable
inductive limit of one-dimensional noncommutative CW complexes, equipped
with a faithful tracial state, and suppose that $K_{0}(E)=\mathbb{Z}[1_{E}]$
and $K_{1}(E)=0$. By \cite[Theorem 3.1]{HirshbergPhillips}, $D*_{\mathrm{red}}E\cong D$.
Combining this with Theorem \ref{thm:freeProductsAbsorb}, every $B=A_{1}*_{\mathrm{red}}A_{2}$
as in that theorem satisfies 
\[
B*_{\mathrm{red}}E\cong B.
\]
Examples of $E$ include $C(T)$ for a finite tree $T$, prime dimension-drop
algebras $Z_{p,q}$ with relatively prime $p,q$, and the Jiang--Su
algebra, with any faithful tracial state. In particular, 
\[
C_{\mathrm{red}}^{*}(\mathbb{F}_{m})*_{\mathrm{red}}E\cong C_{\mathrm{red}}^{*}(\mathbb{F}_{m})\qquad(m\in\{2,3,\ldots,\infty\}).
\]
\end{rem}

\section{Elementary inclusions}\label{sec:elementary-inclusions}

We use the continuous first-order language of tracial von Neumann
algebras; see \cite[Sections 2.3.2--2.5]{FarahHartSherman}. In this
language, formulas are built from traces of $*$-polynomials, continuous
functions of these expressions, and iterated infima and suprema over
norm unit balls. For example, 
\[
\inf_{\|y\|\leq1}\psi(y,a_{1},\ldots,a_{k})
\]
is a formula whenever $\psi$ is one. Every formula is uniformly $L^{2}$-continuous
in its free variables on prescribed operator-norm balls.

An inclusion $N\subset M$ is elementary, written $N\preceq M$, if
\[
\varphi^{N}(a_{1},\ldots,a_{k})=\varphi^{M}(a_{1},\ldots,a_{k})
\]
for every formula $\varphi$ and every tuple from $N$. The superscript
indicates the algebra in which the suprema/infima are evaluated.
\begin{prop}
\label{prop:elementary}If $N\subset M$ has AAP, then the inclusion
is elementary: $N\preceq M$.
\end{prop}

\begin{proof}
By Proposition \ref{prop:aap}, there are surjective trace-preserving
isomorphisms $\Phi_{j}:N\to M$ converging pointwise in $L^{2}$ to
inclusion. Fix a formula $\varphi$ and a tuple $a=(a_{1},\ldots,a_{k})$
from $N$. Isomorphisms preserve formula values: they preserve the
algebraic operations and trace, and surjectivity identifies the domains
of every infimum and supremum. Hence 
\[
\varphi^{N}(a)=\varphi^{M}(\Phi_{j}(a)).
\]
Since $\|\Phi_{j}(a_{i})\|=\|a_{i}\|$ and $\|\Phi_{j}(a_{i})-a_{i}\|_{2}\to0$,
uniform continuity gives 
\[
\varphi^{N}(a)=\lim_{j\to\infty}\varphi^{M}(\Phi_{j}(a))=\varphi^{M}(a).
\]
\end{proof}
In particular, Theorem \ref{thm:free-absorption} and Corollary \ref{cor:countable-absorption}
show that the canonical inclusions 
\[
L(\mathbb{F}_{n})\preceq L(\mathbb{F}_{n+1}),\qquad L(\mathbb{F}_{n})\preceq L(\mathbb{F}_{\infty})\qquad(2\leq n<\infty)
\]
are elementary.

This answers \cite[Question 1.2 and Section 4.2]{GoldbringPi}. Elementarity
of these particular inclusions does not follow from abstract isomorphism
of the free group factors; here it follows from AAP.
\begin{rem}
\label{rem:caap-elementary}The same proof shows that CAAP implies
elementarity in the language of unital $C^{*}$-algebras, using the
norm metric and norm predicates. For trace-preserving CAAP, one may
also name the trace.
\end{rem}

This leaves the question whether every elementary inclusion $N\subset M$
with $N\cong M\cong L(\mathbb{F}_{2})$ must have AAP. There is, however,
a characterization using continuous infinitary logic: for separable
tracial von Neumann algebras, $N\subset M$ has AAP if and only if
it preserves all $L_{\omega_{1}\omega}$-formulas. Here countable
suprema and infima of formulas are allowed subject to a common uniform
continuity modulus. The forward implication follows from the preceding
proof. Conversely, metric Scott analysis \cite[Theorems 3.5 and 5.5]{MetricScottAnalysis}
gives surjective isomorphisms $N\to M$ arbitrarily close to inclusion
on any finite tuple, which yields Proposition \ref{prop:aap}(2).

\appendix

\section{A classical Gaussian rotation and AAP}\label{app:classical-gaussian}

We describe a classical counterpart of our argument; the main difference
is the mechanism for constructing the derivation $\delta$. 

Let $\gamma$ be the standard Gaussian probability measure, with density
$\varphi$ and distribution function $\Phi$. On $([0,1)\times\mathbb{R},dx\,d\gamma(s))$,
write $x,s$ for the two coordinate functions and set
\[
N=L^{\infty}(\sigma(x))\subset M=L^{\infty}(\sigma(x,s)).
\]
Here $x$ is bounded, whereas $s$ is a measurable function belonging
to every $L^{p}$, $1\leq p<\infty$.

\subsection{Constructing the derivation}

For an integer $m\geq1$, put
\[
J_{m}(x)=\lfloor mx\rfloor,\qquad u_{m}=mx-J_{m},\qquad\xi_{m}=\Phi^{-1}(u_{m}).
\]
The variables $J_{m},u_{m}$ are independent, with $J_{m}$ uniform
on $\{0,\ldots,m-1\}$ and $u_{m}$ uniform on $(0,1)$. Thus $J_{m},\xi_{m},s$
are independent, and $\xi_{m},s$ are standard Gaussians.

Using the coordinates $(J_{m},\xi_{m})$, let $\delta_{m}$ be differentiation
in the second coordinate. More precisely, consider functions $f(J_{m},\xi_{m})$
that are smooth and compactly supported in the second variable, and
take the closure of the span of constants and these functions in $L^{2}(N)$.
Gaussian integration by parts on each summand gives
\[
\delta_{m}^{*}1=\xi_{m},\qquad\delta_{m}(\xi_{m})=1,\qquad\delta_{m}(J_{m})=0.
\]
The adjoint here is taken in all of $L^{2}(N)$. The last two identities
hold on the domain of the closure of the derivative operator.

In the original coordinate, this derivation is
\[
\delta_{m}f=a_{m}(x)f'(x),\qquad a_{m}(x)=\frac{\varphi(\xi_{m}(x))}{m}.
\]
Set $a_{m}=0$ at the points $\{k/m:0\leq k\leq m\}$. On each open
interval $[k/m,(k+1)/m]$, $a_{m}'=-\xi_{m}$. Since $x=(J_{m}+\Phi(\xi_{m}))/m$,
we have
\[
\delta_{m}x=a_{m}(x),\qquad\|\delta_{m}x\|_{\infty}\leq\frac{1}{m\sqrt{2\pi}}.
\]
Thus the derivation is small on the bounded generator $x$, while
its conjugate variable $\xi$ is exactly Gaussian and has derivative
$1$.

\subsection{The global flow}

Extend the derivation $\delta_{m}$ to a derivation $D_{m}$ on functions
of the product probability space by
\[
D_{m}f=s\delta_{m}f\quad(f\text{ depending only on }x),\qquad D_{m}s=-\xi_{m}.
\]
In coordinates, this is the single vector field
\[
D_{m}=sa_{m}(x)\partial_{x}-\xi_{m}(x)\partial_{s}=s\partial_{\xi_{m}}-\xi_{m}\partial_{s}.
\]
With $q=\xi_{m}$ and $v=s$, this is the transport operator $B=v\partial_{q}-V'(q)\partial_{v}$
from \cite[Section 7]{VillaniHypocoercivity}, for the quadratic potential
$V(q)=q^{2}/2$. It fixes $J_{m}$, and
\[
D_{m}\xi_{m}=s,\qquad D_{m}s=-\xi_{m}.
\]
In the coordinates $(J_{m},\xi_{m},s)$, its flow is the ordinary
rotation of the two Gaussian coordinates. It is therefore defined
for every real time and preserves their joint law.

Explicitly, away from the null set of mesh points, the flow is
\[
\Psi_{t}^{(m)}(x,s)=\left(\frac{J_{m}(x)+\Phi(\xi_{m}(x)\cos t+s\sin t)}{m},\ s\cos t-\xi_{m}(x)\sin t\right).
\]
Extend it by the identity on that null set. Composition with $\Psi_{t}^{(m)}$
defines a global one-parameter group $\alpha_{t}^{(m)}$ of trace-preserving
automorphisms of $M$. Its extensions to $L^{p}$, $1\leq p<\infty$,
are strongly continuous isometries and satisfy
\[
\alpha_{t}^{(m)}(\xi_{m})=\xi_{m}\cos t+s\sin t,\qquad\alpha_{t}^{(m)}(s)=s\cos t-\xi_{m}\sin t.
\]
These identities concern measurable functions, not bounded elements
of $M$. In particular,
\[
\alpha_{\pi/2}^{(m)}(\xi_{m})=s
\]
holds almost everywhere and in every such $L^{p}$.

There is nevertheless a uniform estimate for the bounded coordinate
$x$. Since $\alpha_{t}^{(m)}$ fixes $J_{m}$,
\[
\alpha_{t}^{(m)}(x)-x=\frac{\Phi(\xi_{m}\cos t+s\sin t)-\Phi(\xi_{m})}{m},
\]
and hence
\[
\|\alpha_{t}^{(m)}(x)-x\|_{\infty}\leq\frac{1}{m}\qquad(t\in\mathbb{R}).
\]
Both values lie in the same interval of length $1/m$.

\subsection{The approximation property}
\begin{prop}
\label{prop:classical-aap}The inclusion $N\subset M$ has AAP.
\end{prop}

\begin{proof}
Take $m=2^{n}$ and write $\beta_{m}=\alpha_{\pi/2}^{(m)}$. Since
$N=L^{\infty}(\sigma(J_{m},\xi_{m}))$, 
\[
\beta_{m}(N)=L^{\infty}(\sigma(J_{m},s)).
\]
For $f\in N$, let $f_{m}=\mathbb{E}[f\mid J_{m}]$. The automorphism
$\beta_{m}$ fixes $f_{m}$, so
\[
\|\beta_{m}(f)-f\|_{2}\leq2\|f-f_{m}\|_{2}\longrightarrow0.
\]
Here we use that the dyadic partitions generate $\sigma(x)$ modulo
null sets.

For $z\in M$, the bounded functions
\[
z_{m}=\mathbb{E}[z\mid J_{m},s]\in\beta_{m}(N)
\]
satisfy $\|z-z_{m}\|_{2}\to0$, because the increasing sigma-algebras
$\sigma(J_{2^{n}},s)$ generate $\sigma(x,s)$. Taking $y_{m}=\beta_{m}^{-1}(z_{m})\in N$
verifies condition (1) of Proposition \ref{prop:aap}.
\end{proof}


\begin{thebibliography}{99}
\bibitem{MetricScottAnalysis}I. Ben Yaacov, M. Doucha, A. Nies and
T. Tsankov, Metric Scott analysis, Adv. Math. 318 (2017), 46-{}-87.

\bibitem{ConnesShlyakhtenko}A. Connes and D. Shlyakhtenko, $L^{2}$-homology
for von Neumann algebras, J. Reine Angew. Math. 586 (2005), 125-{}-168.

\bibitem{FarahHartSherman}I. Farah, B. Hart and D. Sherman, Model
theory of operator algebras II: Model theory, Israel J. Math. 201
(2014), 477-{}-505.

\bibitem{GaboriauL2}D. Gaboriau, Invariants $\ell^{2}$ de relations
d’équivalence et de groupes, Publ. Math. IHÉS 95 (2002), 93--150.

\bibitem{GaoKunnawalkamElayavalli}D. Gao and S. Kunnawalkam Elayavalli,
Non-isomorphism of reduced free group $C^{*}$-algebras, arXiv:2602.09995
(2026).

\bibitem{GoldbringPi}I. Goldbring and J. Pi, On the first-order free
group factor elementary equivalence, J. Operator Theory 94 (2025),
no. 1, 129-{}-150.

\bibitem{HirshbergPhillips}I. Hirshberg and N. C. Phillips, An isomorphism
theorem for infinite reduced free products, arXiv:2602.10220v2 (2026).

\bibitem{NicaSpeicher}A. Nica and R. Speicher, Lectures on the Combinatorics
of Free Probability, London Mathematical Society Lecture Note Series
335, Cambridge University Press, 2006.

\bibitem{OpenAIFreeGroupFactors}OpenAI, An isomorphism of the free
group factors, preprint, September 23, 2026, family 287 in the collection
at github.com/openai/math.

\bibitem{PimsnerVoiculescu}M. Pimsner and D. Voiculescu, K-groups
of reduced crossed products by free groups, J. Operator Theory 8 (1982),
no. 1, 131-{}-156.

\bibitem{PopaShlyakhtenkoGt}S. Popa and D. Shlyakhtenko, Representing
interpolated free group factors as group factors, Groups Geom. Dyn.
14 (2020), no. 3, 837--855.

\bibitem{ShlyakhtenkoMicrostates}D. Shlyakhtenko, Lower estimates
on microstates free entropy dimension, Anal. PDE 2 (2009), no. 2,
119-{}-146.

\bibitem{ShlyakhtenkoFuchsian}D. Shlyakhtenko, On the II$_{1}$ Factors
of Fuchsian Groups, arXiv:2609.11074v1 (2026).

\bibitem{ShlyakhtenkoPartialConjugates}D. Shlyakhtenko, Some estimates
for non-microstates free entropy dimension with applications to $q$-semicircular
families, Int. Math. Res. Not. 2004 (2004), no. 51, 2757-{}-2772.

\bibitem{ShlyakhtenkoSofic}D. Shlyakhtenko, Von Neumann algebras
of sofic groups with $\beta_{1}^{(2)}=0$ are strongly $1$-bounded.
Journal of Operator Theory 85 (2021), no. 1, 217--228.

\bibitem{ThielHaarUnitaries}H. Thiel, Diffuse traces and Haar unitaries,
Amer. J. Math. 146 (2024), no. 5, 1305-{}-1337.

\bibitem{VillaniHypocoercivity}C. Villani, Hypocoercivity, Mem. Amer.
Math. Soc. 202 (2009), no. 950, iv+141 pp.; arXiv:math/0609050.

\bibitem{VoiculescuEntropyII}D. Voiculescu, The analogues of entropy
and of Fisher's information measure in free probability theory, II.
Invent. Math. 118 (1994), 411--440.

\bibitem{VoiculescuEntropyV}D. Voiculescu, The analogues of entropy
and of Fisher's information measure in free probability theory. V.
Noncommutative Hilbert transforms, Invent. Math. 132 (1998), no. 1,
189-{}-227.

\bibitem{VoiculescuEntropyVI}D. Voiculescu, The analogues of entropy
and of Fisher's information measure in free probability theory. VI.
Liberation and mutual free information, Adv. Math. 146 (1999), no.
2, 101-{}-166.

\bibitem{VoiculescuCyclomorphy}D. Voiculescu, Cyclomorphy, Int. Math.
Res. Not. 2002 (2002), no. 6, 299-{}-332.

\bibitem{VoiculescuSymmetries}D. Voiculescu, “Symmetries of some
reduced free product $C^{*}$-algebras,” in Operator Algebras and
Their Connections with Topology and Ergodic Theory (Buşteni, 1983),
Lecture Notes in Mathematics 1132, Springer, Berlin, 1985, pp. 556--588.

\end{thebibliography}
\end{document}